\documentclass[11pt]{article}
\usepackage{amsfonts}
\usepackage{amsfonts}
\usepackage{amsfonts}
\usepackage{mathrsfs}
\usepackage{bm}
\usepackage{cite}
\usepackage[colorlinks,linkcolor=blue,citecolor=blue]{hyperref}
\usepackage{amssymb,amsmath} 
 \allowdisplaybreaks

\catcode`!=11
\let\!int\int \def\int{\displaystyle\!int}
\let\!lim\lim \def\lim{\displaystyle\!lim}
\let\!sum\sum \def\sum{\displaystyle\!sum}
\let\!sup\sup \def\sup{\displaystyle\!sup}
\let\!inf\inf \def\inf{\displaystyle\!inf}
\let\!cap\cap \def\cap{\displaystyle\!cap}
\let\!max\max \def\max{\displaystyle\!max}
\let\!min\min \def\min{\displaystyle\!min}
\let\!frac\frac \def\frac{\displaystyle\!frac}
\catcode`!=12

\let\oldsection\section
\renewcommand\section{\setcounter{equation}{0}\oldsection}

\allowdisplaybreaks
\def\pf{\it{Proof.}\rm\quad}

\def\N{\mathbb{N}}

\newtheorem{thm}{Theorem}[section]
\newtheorem{lem}[thm]{Lemma}
\newtheorem{cor}[thm]{Corollary}
\newtheorem{con}[thm]{Conjecture}
\newtheorem{re}{Remark}[section]

\begin{document}
\title {\bf Evaluations of nonlinear Euler sums of weight ten}
\author{
{Ce Xu\thanks{Corresponding author. Email: 15959259051@163.com}}\\[1mm]
\small School of Mathematical Sciences, Xiamen University\\
\small Xiamen
361005, P.R. China\\
}

\date{}
\maketitle \noindent{\bf Abstract } In this paper we present a new family of identities for Euler sums and integrals of polylogarithms by using the methods of generating function and integral representations of series. Then we apply it to obtain the closed forms of all quadratic Euler sums of weight equal to ten. Furthermore, we also establish some relations between multiple zeta (star) values and nonlinear Euler sums. As applications of these relations, we give new closed form representations of several cubic Euler sums through single zeta values and linear sums. Finally, with the help of numerical computations of Mathematica or Maple, we evaluate several other Euler sums of weight ten.
\\[2mm]
\\[2mm]
\tableofcontents
\section{Introduction}
Let $k$ and $n$ be positive integers. The generalized harmonic numbers are defined by the finite sums
\[
H^{(k)}_n:=\sum\limits_{j=1}^n\frac {1}{j^k}\quad {\rm and}\quad H^{(k)}_0:=0.
\]
When $k=1$, then we let $H_n=H^{(1)}_n$, which are called the classical harmonic numbers. Letting $n$ tend to infinity in above definition, then the generalized harmonic number $H^{(k)}_n$ converges to the Riemann zeta value $\zeta(k)$, i.e., the values of the Riemann zeta function $\zeta (s)\ (\Re(s)>1)$ at the positive integer arguments:
 \[\mathop {\lim }\limits_{n \to \infty } H_n^{\left( k \right)} = \zeta( k),\ {\Re} \left( k \right) > 1,\ k\in \N.\]
The classical Euler sum $S_{p,q}$ (also called the linear Euler sum or nested sum of zeta value) is defined by the convergent series (see Berndt \cite{B1985,B1989})
\[S_{p,q}:=\sum\limits_{n = 1}^\infty  {\frac{{H_n^{\left( p \right)}}}
{{{n^q}}}}, \]
where $p,q$ are positive integers with $q \geq 2$, and $w:=p+q$ denotes the weight of linear sums $S_{p,q}$.
Flajolet and Salvy used the harmonic numbers $H^{(k)}_n$ to define the generalized Euler sums ${S_{{\bf S},q}}$, which are defined by (\cite{FS1998})
\[{S_{{\bf S},q}} := \sum\limits_{n = 1}^\infty  {\frac{{H_n^{\left( {{s_1}} \right)}H_n^{\left( {{s_2}} \right)} \cdots H_n^{\left( {{s_r}} \right)}}}
{{{n^q}}}}.\]
Here ${\bf S}:=(s_1,s_2,\ldots,s_r)\quad(r,s_i\in \N, i=1,2,\ldots,r)$ with $s_1\leq s_2\leq \ldots\leq s_r$ and $q\geq 2$. The quantity $w:={s _1} +  \cdots  + {s _r} + q$ is called the weight and the quantity $r$ is called the degree of the sum.  Since repeated summands in partitions are indicated by powers, we denote, for example, the sum
 \[{S_{{1^2}{2^3}4,q}} = {S_{112224,q}} = \sum\limits_{n = 1}^\infty  {\frac{{H_n^2[H^{(2)} _n]^3{H^{(4)} _n}}}{{{n^q}}}}. \]
When $r>1$, then we call the sums ${S_{{\bf S},q}}$ the nonlinear Euler sums. In particular, they are called the quadratic and cubic Euler sums when $r=2$ and $3$, respectively.

Investigation of Euler sums has a long history. Euler's original contribution was a method to reduce linear sums to certain rational linear combinations of products of zeta values. For example, Euler proved that the linear sums $S_{p,q}$ are reducible to zeta values whenever $p+q$ is less than 7 or when $p+q$ is odd and less than 13. Examples for such evaluations, all due to Euler, are
\begin{align*}
&{S_{1,4}} = 3\zeta \left( 5 \right) - \zeta \left( 2 \right)\zeta \left( 3 \right),\\
&{S_{2,3}} =  - \frac{9}{2}\zeta \left( 5 \right) + 3\zeta \left( 2 \right)\zeta \left( 3 \right),\\
&{S_{4,2}} = \frac{{37}}{{12}}\zeta \left( 6 \right) - {\zeta ^2}\left( 3 \right),\\
&{S_{3,4}} = 18\zeta \left( 7 \right) - 10\zeta \left( 2 \right)\zeta \left( 5 \right),\\
&{S_{2,5}} = 5\zeta \left( 2 \right)\zeta \left( 5 \right) + 2\zeta \left( 3 \right)\zeta \left( 4 \right) - 10\zeta \left( 7 \right).
\end{align*}
The linear sums $S_{p,q}$ can be evaluated in terms of zeta values in the following cases: $p=1,p=q,p+q$ odd and $p+q=6$ with $q\geq 2$\ (for more details, see \cite{BBG1994,BBG1995,FS1998}). For example, a evaluation of Euler sums $S_{p,q}$ of odd weight $p+q$ in terms of zeta values can be accomplished through Theorem 3.1 of \cite{FS1998}.  At weight $w=8$, it appears that $S_{2,6}$ can not be reduced to zeta values and their products, though we have no way of proving that such a reduction cannot exist. We cannot even prove that $\zeta(3)/\pi^3$ is irrational. I shall take $S_{2,6}$ as an irreducible linear sums of weight 8. Then all other linear sums of weight 8 may be reduced to $S_{2,6}$ and zeta values. For example,
\[5{S_{2,6}} + 2{S_{3,5}} =  - \frac{{21}}{4}\zeta \left( 8 \right) + 10\zeta \left( 3 \right)\zeta \left( 5 \right).\]
It is proven that the number of irreducible linear sums $S_{p,q}$ of even weight $w=2n$  is no greater than $\left\lceil {n/3} \right\rceil  - 1$. Upon to weight$=10$,  we may take the irreducible linear sums to be $S_{2,6}$ and $S_{2,8}$. For instance, we have the following relations
\begin{align*}
&7{S_{2,8}} + 2{S_{3,7}} =  - \frac{{33}}{2}\zeta \left( {10} \right) + 14\zeta \left( 3 \right)\zeta \left( 7 \right) + 8{\zeta ^2}\left( 5 \right),\\
&7{S_{2,8}} - 2{S_{4,6}} =  - \frac{{227}}{{10}}\zeta \left( {10} \right) + 14\zeta \left( 3 \right)\zeta \left( 7 \right) + 10{\zeta ^2}\left( 5 \right).
\end{align*}

Similarly, many values of nonlinear Euler sums $S_{{\bf S},q}$ also can be expressed as a rational linear combination of zeta values and linear sums. For example, in 1994, Bailey et al. \cite{BBG1994} proved that all Euler sums of the form $S_{1^p,q}$ for weights $p+q\in \{3,4,5,6,7,9\}$ are reducible to {\bf Q}-linear combinations of zeta values by using the experimental method. In 1995, Borwein et al. \cite{BBG1995} showed that the quadratic sums $S_{ 1^2,q}$ can reduce to linear sums $S_{2,q}$ and polynomials in zeta values. In 1998, Flajolet and Salvy \cite{FS1998} used the contour integral representations and residue computation to show that the quadratic sums $S_{p_1p_2,q}$ are reducible to linear sums and zeta values when the weight $p_1 + p_2 + q$ is even and $p_1,p_2>1$. Moreover, in \cite{X2017,Xu2017}, we proved that all Euler sums of weight $\leq 8$ are reducible to $\mathbb{Q}$-linear combinations of single zeta monomials with the addition of $\{S_{2,6}\}$ for weight 8. For weight 9, all Euler sums of the form ${S_{{s_1} \cdots {s_k},q}}$ with $q\in \{4,5,6,7\}$ are expressible polynomially in terms of zeta values. For weight $1+p+q=10$, all quadratic sums $S_{1p_,q}$ are reducible to $S_{2,6}$ and $S_{2,8}$. Very recently, Wang et al \cite{W2017} shown that all Euler sums of weight eight are reducible to linear sums, and proved that all Euler sums of weight nine are reducible to zeta values. Some simple examples are as follows
\begin{align*}
&S_{1^3,4} =
\frac{{231}}{{16}}\zeta (7) - \frac{{51}}{4}\zeta (3)\zeta (4) +
2\zeta (2)\zeta (5),\\
&S_{2^2,3}  =  - \frac{{155}}{8}\zeta(7)
+ \frac{{19}}{2}\zeta(3)\zeta( 4) + 5\zeta( 2 )\zeta( 5),\\
&{S_{{1^4}2,2}} =
- \frac{{6631}}{{288}}\zeta \left( 8 \right) + 90\zeta \left( 3 \right)\zeta \left( 5 \right) + 3\zeta \left( 2 \right){\zeta ^2}\left( 3 \right) - \frac{{47}}{2}{S_{2,6}},\\
&{S_{123,2}}
=  - \frac{{181}}{{288}}\zeta \left( 8 \right) + \frac{{15}}{2}\zeta \left( 3 \right)\zeta \left( 5 \right) - \frac{3}{2}\zeta \left( 2 \right){\zeta ^2}\left( 3 \right) - \frac{7}{4}{S_{2,6}},\\
&S_{2^2,5} =
 - \frac{{1069}}{{36}}\zeta( 9) + \frac{4}{3}{\zeta ^3}( 3) + 7\zeta( 2)\zeta( 7)
 - \frac{4}{3}\zeta ( 3)\zeta( 6) + \frac{{33}}{2}\zeta( 4)\zeta( 5),\\
&S_{25,2}
 = \frac{{2059}}{{72}}\zeta( 9) - \frac{2}{3}{\zeta ^3}( 3) - 14\zeta( 2)\zeta( 7)
 + \frac{8}{3}\zeta( 3 )\zeta( 6) - 5\zeta( 4)\zeta( 5),\\
&S_{17,2}= \frac{{242}}
{{15}}\zeta( {10}) - \frac{{25}}
{4}{\zeta ^2}( 5) - \frac{{19}}
{2}\zeta( 3)\zeta( 7)+{\zeta ^2}( 3)\zeta( 4)+\zeta( 2){S_{2,6}} + \frac{5}
{4}{S_{2,8}}.
\end{align*}
The relationship between the values of the Riemann zeta function and Euler sums has been
studied by many authors. Besides the works referred to above, there are many other researches devoted to the Euler sums. For details and historical introductions, please see \cite{BBG1994,BBG1995,BZB2008,FS1998,F2005,M2014,Xu2016,X2016,X2017,Xyz2016,Xu2017} and references therein.

Since the subject of this paper are Euler sums and multiple zeta (star) values. Hence, we now give a brief introductions of the multiple harmonic (star) sums and the multiple zeta (star) values. Let $s_1,\ldots,s_k$ be positive integers. The multiple harmonic star sums (MHSS) and multiple harmonic sums (MHS) are defined by (\cite{KP2013,Xu2017})
\[{\zeta _n}\left( {{s_1},{s_2}, \cdots ,{s_k}} \right): = \sum\limits_{n \ge {n_1} > {n_2} >  \cdots  > {n_k} \ge 1} {\frac{1}{{n_1^{{s_1}}n_2^{{s_2}} \cdots n_k^{{s_k}}}}} ,\]
\[{\zeta_n ^ \star }\left( {{s_1},{s_2}, \cdots ,{s_k}} \right): = \sum\limits_{n \ge {n_1} \ge {n_2} \ge  \cdots  \ge {n_k} \ge 1} {\frac{1}{{n_1^{{s_1}}n_2^{{s_2}} \cdots n_k^{{s_k}}}}},\]
when $n<k$, then ${\zeta _n}\left( {{s_1},{s_2}, \cdots ,{s_k}} \right)=0$, and ${\zeta _n}\left(\emptyset \right)={\zeta^\star _n}\left(\emptyset \right)=1$. The integers $k$ and $w:=s_1+\ldots+s_k$ are called the depth and the weight of a multiple harmonic (star) sum. For convenience, by ${\left\{ {{s_1}, \ldots ,{s_j}} \right\}_d}$ we denote the sequence of depth $dj$ with $d$ repetitions of ${\left\{ {{s_1}, \ldots ,{s_j}} \right\}}$. For example,
\[\zeta_n\left( {5,3,{{\left\{ 1 \right\}}_2}} \right) = \zeta_n\left( {5,3,1,1} \right),\;{\zeta_n^ \star }\left( {4,2,{{\left\{ 1 \right\}}_3}} \right) = {\zeta_n ^ \star }\left( {4,2,1,1,1} \right).\]
 When taking the limit $n\rightarrow \infty$ we get the so-called the multiple zeta value (MZV) and the multiple zeta star value (MZSV), respectively (\cite{BBV2010,BBBL1997,Bro2013,DZ2012,Z2015}):
\begin{align*}
&\zeta \left( {{s_1},{s_2}, \cdots ,{s_k}} \right) = \mathop {\lim }\limits_{n \to \infty } \zeta \left( {{s_1},{s_2}, \cdots ,{s_k}} \right),\\
&\zeta _{}^ \star \left( {{s_1},{s_2}, \cdots ,{s_k}} \right) = \mathop {\lim }\limits_{n \to \infty } \zeta _{}^ \star \left( {{s_1},{s_2}, \cdots ,{s_k}} \right)
\end{align*}
defined for $s_2,\ldots,s_k\geq 1$ and $s_1\geq 2$ to ensure convergence of the series.

Similarly, the alternating multiple harmonic (star) sums are closely related to the MHSS and MHS, which are defined by
\begin{align*}
&\zeta_n \left( \bf s \right)\equiv\zeta_n \left( {{s_1}, \ldots ,{s_k}} \right): = \sum\limits_{n\geq {n_1} >  \cdots  > {n_k} > 0} {\prod\limits_{j = 1}^k {n_j^{ - \left| {{s_j}} \right|}}{\rm sgn}(s_j)^{n_j}} ,\\
&\zeta_n^\star \left( \bf s \right)\equiv{\zeta_n ^ \star }\left( {{s_1}, \ldots ,{s_k}} \right): = \sum\limits_{n\geq{n_1} \ge  \cdots  \ge {n_k} \ge 1} {\prod\limits_{j = 1}^k {n_j^{ - \left| {{s_j}} \right|}{\rm sgn}(s_j)^{n_j}} } ,
\end{align*}
where $s_j$ stands for non-zero integer, and
\[{\mathop{\rm sgn}} \left( {{s_j}} \right): = \left\{ {\begin{array}{*{20}{c}}
   {1,} & {{s_j} > 0,}  \\
   { - 1,} & {{s_j} < 0.}  \\
\end{array}} \right.\]
 We may compactly indicate the presence of an alternating sign. When ${\rm sgn}(s_j)=-1$,  by placing a bar over the
corresponding integer exponent $s_j$. Thus we write
\[\zeta _n^ \star \left( {\bar 2,3} \right) =\zeta _n^ \star \left( {-2,3} \right)= \sum\limits_{n \ge {n_1} \ge {n_2} \ge 1}^{} {\frac{{{{\left( { - 1} \right)}^{{n_1}}}}}{{n_1^2{n_2^3}}}},\]
\[{\zeta _n}\left( {\bar 2,3,\bar 1,4} \right)={\zeta _n}\left( {-2,3,- 1,4} \right) = \sum\limits_{n \ge {n_1} > {n_2} > {n_3} > {n_4} \ge 1}^{} {\frac{{{{\left( { - 1} \right)}^{{n_1} + {n_3}}}}}{{n_1^2n_2^3{n_3}n_4^4}}} .\]
Clearly, the limit cases of alternating multiple harmonic (star) sums give rise to alternating multiple zeta (star) values, for example
\begin{align*}
&{\zeta ^ \star }\left( {\bar 2,3,5} \right) = \mathop {\lim }\limits_{n \to \infty } \zeta _n^ \star \left( {\bar 2,3,5} \right) = \sum\limits_{{n_1} \ge {n_2} \ge 1}^{} {\frac{{{{\left( { - 1} \right)}^{{n_1}}}}}{{n_1^2n_2^3n_3^5}}} \\
&\zeta \left( {\bar 2,3,\bar 1,4} \right) = \mathop {\lim }\limits_{n \to \infty } {\zeta _n}\left( {\bar 2,3,\bar 1,4} \right) = \sum\limits_{{n_1} > {n_2} > {n_3} > {n_4} \ge 1}^{} {\frac{{{{\left( { - 1} \right)}^{{n_1} + {n_3}}}}}{{n_1^2n_2^3{n_3}n_4^4}}} .
\end{align*}
In \cite{Xu2017}, the author proved the following relations
\[{\zeta _n}\left( {{{\left\{ 1 \right\}}_k}} \right) = \frac{{s\left( {n + 1,k + 1} \right)}}
{{n!}},\;\zeta _n^ \star \left( {{{\{ 1\} }_k}} \right) = \frac{{{Y_k}\left( n \right)}}
{{k!}},\;\left( {n,k \in \N} \right)\]
where $s(n,k)$ and $Y_k(n):= {Y_k}\left( {{H _n}( 1 ),1!{H^{(2)} _n},2!{H^{(3)} _n}, \cdots \left( {k - 1} \right)!{H^{(k)} _n}} \right)$ denote the (unsigned) Stirling number of the first kind and Bell number, respectively. Here ${Y_k}\left( {{x_1},{x_2}, \cdots } \right)$ stands for the complete exponential Bell polynomial is defined by (see \cite{L1974})
\[\exp \left( {\sum\limits_{m \ge 1}^{} {{x_m}\frac{{{t^m}}}{{m!}}} } \right) = 1 + \sum\limits_{k \ge 1}^{} {{Y_k}\left( {{x_1},{x_2}, \cdots } \right)\frac{{{t^k}}}{{k!}}}.\]
 From the definition of the complete exponential Bell polynomial, we deduce that
$${Y_1}\left( n \right) = {H_n},{Y_2}\left( n \right) = H_n^2 + {H^{(2)} _n},{Y_3}\left( n \right) =  H_n^3+ 3{H_n}{H^{(2)} _n}+ 2{H^{(3)} _n},$$
\[{Y_4}\left( n \right) = H_n^4 + 8{H_n}{H^{(3)} _n} + 6H_n^2{H^{(2)} _n} + 3(H^{(2)} _n)^2 + 6{H^{(4)} _n},\]
\[{Y_5}\left( n \right) = H_n^5 + 10H_n^3{H^{(2)} _n} + 20H_n^2{H^{(3)}_n} + 15{H_n}({H^{(2)}_n})^2 + 30{H_n}{H^{(4)} _n}+ 20{H^{(2)} _n}{H^{(3)} _n} + 24{H^{(5)} _n},\]

Similarly, using the a bar notation, we also define the generalized alternating Euler sums ${{\bar S}_{{\bf S},q}}$ by
\[{{\bar S}_{{\bf{S}},q}}: = \sum\limits_{n = 1}^\infty  {\frac{{H_n^{\left( {{s_1}} \right)}H_n^{\left( {{s_2}} \right)} \cdots H_n^{\left( {{s_r}} \right)}}}{{{n^q}}}{{\left( { - 1} \right)}^{n - 1}}},\]
Here ${\bf S}:=(s_1,s_2,\ldots,s_r)\quad(r,s_i\in \N, i=1,2,\ldots,r)$ with $s_1\leq s_2\leq \ldots\leq s_r$ and $q\geq 1$.
By using the above notations, the following relations between linear Euler sums and multiple zeta star values are obviously
\begin{align*}
&{S_{p,q}} = {\zeta ^ \star }\left( {q,p} \right) = \zeta \left( {q,p} \right) + \zeta \left( {q + p} \right),\\
&{{\bar S}_{p,q}} =  - {\zeta ^ \star }\left( {\bar q,p} \right) =  - \zeta \left( {\bar q,p} \right) - \zeta \left( {\overline {q + p}} \right).
\end{align*}
Note that a fast program for evaluating MZVs has been developed at the CECM (Centre for Experimental and Constructive Mathematics, Simon Fraser University)  and is available for public use via the World Wide Web interface called ¡°EZ Face¡± (an abbreviation for Euler Zetas interFace) at the URL
\url{http://wayback.cecm.sfu.ca/projects/EZFace/}.
This publicly accessible interface currently allows one to evaluate the sums
\[{\rm z}\left( {{s_1}, \ldots ,{s_k}} \right): = \sum\limits_{{n_1} >  \cdots  > {n_k} > 0} {\prod\limits_{j = 1}^k {n_j^{ - \left| {{s_j}} \right|}} {\rm{sgn}}{{({s_j})}^{{n_j}}}} \]
for non-zero integers $s_1,\ldots,s_k$, and
\[{\rm zp}\left( {p,{s_1}, \ldots ,{s_k}} \right): = \sum\limits_{{n_1} >  \cdots  > {n_k} > 0} {{p^{ - {n_1}}}\prod\limits_{j = 1}^k {n_j^{ - {s_j}}} } \]
for real $p\geq1$ and positive integers $s_1,\ldots,s_k$.

A good deal of work on multiple zeta (star) values has focused on the problem of determining when `complicated' sums can be expressed in terms of `simpler' sums. Thus, researchers are interested in determining which sums can be expressed in terms of other sums of lesser depth. The MZVs of depth=3 are investigated \cite{BG1996,FS1998}. In \cite{BG1996}, Borwein and Girgensohn proved that all $\zeta \left( {s_1,s_2,s_3} \right)$ with $s_1+s_2+s_3$ is even or less than or equal to 10 or $s_1+s_2+s_3=12$ were reducible to zeta values and double sums. There are many works investigating multiple zeta values, see for example \cite{BBV2010,BBBL1997,Bro2013,EW2012,H1992,KP2013,DZ2012,Z2015}, and references therein.

In this paper, we will establish some explicit relations involving Euler sums, polylogarithms' integrals and multiple zeta (star) values. Then, we use these relations to evaluate several nonlinear Euler sums of weight=10. Finally, we give all numerical values of Euler sums with weight=10. The obtained results agree well with the simulation data, see Table 1 and 2.
\section{Main Results}
 We will use our equation system to obtain the following explicit evaluations of nonlinear Euler sums of weight equal to ten.
\begin{align}\label{equ:1}
{S_{26,2}} =&  - \frac{{2997}}{{80}}\zeta \left( {10} \right) + 23{\zeta ^2}\left( 5 \right) + 35\zeta \left( 3 \right)\zeta \left( 7 \right) - 8\zeta \left( 2 \right)\zeta \left( 3 \right)\zeta \left( 5 \right)\nonumber\\
& - 2{\zeta ^2}\left( 3 \right)\zeta \left( 4 \right) - \frac{{13}}{2}{S_{2,8}} - \zeta \left( 2 \right){S_{2,6}},
\end{align}
\begin{align}\label{equ:2}
{S_{23,5}}  = & - \frac{{2227}}{{32}}\zeta \left( {10} \right) + \frac{{89}}{2}{\zeta ^2}\left( 5 \right) + 56\zeta \left( 3 \right)\zeta \left( 7 \right) - 15\zeta \left( 2 \right)\zeta \left( 3 \right)\zeta \left( 5 \right)\nonumber\\
& - 10{S_{2,8}} - \frac{5}{2}\zeta \left( 2 \right){S_{2,6}},
\end{align}
\begin{align}\label{equ:3}
{S_{25,3}}  =& \frac{{3223}}{{160}}\zeta \left( {10} \right) - \frac{{17}}{4}{\zeta ^2}\left( 5 \right) + \frac{{63}}{2}\zeta \left( 3 \right)\zeta \left( 7 \right) - 32\zeta \left( 2 \right)\zeta \left( 3 \right)\zeta \left( 5 \right)\nonumber\\
& - 2{\zeta ^2}\left( 3 \right)\zeta \left( 4 \right) - \frac{{19}}{4}{S_{2,8}}{\rm{ + }}\frac{{25}}{2}\zeta \left( 2 \right){S_{2,6}},
\end{align}
\begin{align}\label{equ:4}
{S_{35,2}}  = & \frac{{407}}{{10}}\zeta \left( {10} \right) - \frac{{143}}{4}{\zeta ^2}\left( 5 \right) - \frac{{161}}{2}\zeta \left( 3 \right)\zeta \left( 7 \right) + 48\zeta \left( 2 \right)\zeta \left( 3 \right)\zeta \left( 5 \right)\nonumber\\
& + 2{\zeta ^2}\left( 3 \right)\zeta \left( 4 \right) + \frac{{49}}{4}{S_{2,8}} - 10\zeta \left( 2 \right){S_{2,6}},
\end{align}
\begin{align}\label{equ:5}
{S_{24,4}}   =& \frac{{1519}}{{120}}\zeta \left( {10} \right) - 20{\zeta ^2}\left( 5 \right) - 56\zeta \left( 3 \right)\zeta \left( 7 \right) + 40\zeta \left( 2 \right)\zeta \left( 3 \right)\zeta \left( 5 \right)\nonumber\\
& + {\zeta ^2}\left( 3 \right)\zeta \left( 4 \right) + 11{S_{2,8}} - 10\zeta \left( 2 \right){S_{2,6}},
\end{align}
\begin{align}\label{equ:6}
{S_{{4^2},2}} = & - \frac{{203}}{{120}}\zeta \left( {10} \right) + 30{\zeta ^2}\left( 5 \right) + 98\zeta \left( 3 \right)\zeta \left( 7 \right) - 80\zeta \left( 2 \right)\zeta \left( 3 \right)\zeta \left( 5 \right)\nonumber\\
&- 2{\zeta ^2}\left( 3 \right)\zeta \left( 4 \right) - 14{S_{2,8}} + 20\zeta \left( 2 \right){S_{2,6}},
\end{align}
\begin{align}\label{equ:7}
{S_{{3^2},4}}   = & \frac{{939}}{{10}}\zeta \left( {10} \right) - 46{\zeta ^2}\left( 5 \right) + 7\zeta \left( 3 \right)\zeta \left( 7 \right) - 40\zeta \left( 2 \right)\zeta \left( 3 \right)\zeta \left( 5 \right)\nonumber\\
& - \frac{7}{2}{S_{2,8}} + 20\zeta \left( 2 \right){S_{2,6}},
\end{align}
\begin{align}\label{equ:8}
{S_{34,3}}  = & - \frac{{2001}}{{40}}\zeta \left( {10} \right) + \frac{{49}}{2}{\zeta ^2}\left( 5 \right) + 20\zeta \left( 2 \right)\zeta \left( 3 \right)\zeta \left( 5 \right)\nonumber\\
&+ \frac{1}{2}{\zeta ^2}\left( 3 \right)\zeta \left( 4 \right) - 10\zeta \left( 2 \right){S_{2,6}},
\end{align}
\begin{align}\label{equ:9}
{S_{{2^3},4}}  =& \frac{{458433}}{{640}}\zeta \left( {10} \right) - 190{\zeta ^2}\left( 5 \right) - 374\zeta \left( 3 \right)\zeta \left( 7 \right) + 24\zeta \left( 2 \right)\zeta \left( 3 \right)\zeta \left( 5 \right)\nonumber\\
& + \frac{{57}}{2}{\zeta ^2}\left( 3 \right)\zeta \left( 4 \right) + 3\zeta \left( 2 \right){S_{2,6}} + 58{S_{2,8}} - 24{{\bar S}_{2,8}} - 192{{\bar S}_{1, 9}},
\end{align}
\begin{align}\label{equ:10}
{S_{{2^4},2}}  = & \frac{{415987}}{{320}}\zeta \left( {10} \right) - 298{\zeta ^2}\left( 5 \right) - 622\zeta \left( 3 \right)\zeta \left( 7 \right) + 16\zeta \left( 2 \right)\zeta \left( 3 \right)\zeta \left( 5 \right)\nonumber\\
& + 49{\zeta ^2}\left( 3 \right)\zeta \left( 4 \right) + 2\zeta \left( 2 \right){S_{2,6}} + 94{S_{2,8}} - 48{{\bar S}_{2, 8}} - 384{{\bar S}_{1,9}},
\end{align}
\begin{align}\label{equ:11}
{S_{{2^2}4,2}} =&  - \frac{{102739}}{{640}}\zeta \left( {10} \right) + 4{\zeta ^2}\left( 5 \right) + 8\zeta \left( 3 \right)\zeta \left( 7 \right) + 48\zeta \left( 2 \right)\zeta \left( 3 \right)\zeta \left( 5 \right)\nonumber\\
& - \frac{9}{2}{\zeta ^2}\left( 3 \right)\zeta \left( 4 \right) - 9\zeta \left( 2 \right){S_{2,6}} + {S_{2,8}} + 8{{\bar S}_{2, 8}} + 64{{\bar S}_{1, 9}},
\end{align}
\begin{align}\label{equ:12}
 {S_{{1^4},6}}  =& \frac{{5657}}{{40}}\zeta \left( {10} \right) - \frac{{463}}{8}{\zeta ^2}\left( 5 \right) - \frac{{765}}{8}\zeta \left( 3 \right)\zeta \left( 7 \right) - 12\zeta \left( 2 \right)\zeta \left( 3 \right)\zeta \left( 5 \right)\nonumber\\
 & + \frac{{51}}{2}{\zeta ^2}\left( 3 \right)\zeta \left( 4 \right) + 2\zeta \left( 2 \right){S_{2,6}} + \frac{{303}}{{16}}{S_{2,8}},
\end{align}
\begin{align}\label{equ:13}
{S_{{1^2}2,6}} = & - \frac{{5107}}{{320}}\zeta \left( {10} \right) + \frac{{147}}{{16}}{\zeta ^2}\left( 5 \right) + \frac{{345}}{{16}}\zeta \left( 3 \right)\zeta \left( 7 \right) - 6\zeta \left( 2 \right)\zeta \left( 3 \right)\zeta \left( 5 \right)\nonumber\\
& - \frac{{13}}{4}{\zeta ^2}\left( 3 \right)\zeta \left( 4 \right) - \frac{{51}}{{32}}{S_{2,8}},
\end{align}
\begin{align}\label{equ:14}
{S_{{1^2}3,5}} = & \frac{{213}}{{160}}\zeta \left( {10} \right) - \frac{{43}}{8}{\zeta ^2}\left( 5 \right) - \frac{{97}}{8}\zeta \left( 3 \right)\zeta \left( 7 \right) + 13\zeta \left( 2 \right)\zeta \left( 3 \right)\zeta \left( 5 \right)\nonumber\\
& - \frac{1}{4}{\zeta ^2}\left( 3 \right)\zeta \left( 4 \right) - \zeta \left( 2 \right){S_{2,6}} - \frac{{69}}{{16}}{S_{2,8}},
\end{align}
\begin{align}\label{equ:15}
{S_{{1^2}4,4}}   = & \frac{{7749}}{{160}}\zeta \left( {10} \right) - 16{\zeta ^2}\left( 5 \right) - \frac{{125}}{8}\zeta \left( 3 \right)\zeta \left( 7 \right) - 14\zeta \left( 2 \right)\zeta \left( 3 \right)\zeta \left( 5 \right)\nonumber\\
&+ \frac{3}{2}{\zeta ^2}\left( 3 \right)\zeta \left( 4 \right) + 11{S_{2,8}} + \frac{5}{2}\zeta \left( 2 \right){S_{2,6}},
\end{align}
\begin{align}\label{equ:16}
{S_{{1^2}5,3}}   = & - \frac{{2675}}{{96}}\zeta \left( {10} \right) + \frac{{133}}{8}{\zeta ^2}\left( 5 \right) + \frac{{71}}{8}\zeta \left( 3 \right)\zeta \left( 7 \right) + 5\zeta \left( 2 \right)\zeta \left( 3 \right)\zeta \left( 5 \right)\nonumber\\
& + \frac{3}{4}{\zeta ^2}\left( 3 \right)\zeta \left( 4 \right) - \frac{{167}}{{16}}{S_{2,8}},
\end{align}
\begin{align}\label{equ:17}
{S_{{1^2}6,2}}  =& \frac{{7193}}{{320}}\zeta \left( {10} \right) - \frac{{179}}{{16}}{\zeta ^2}\left( 5 \right) - \frac{{275}}{{16}}\zeta \left( 3 \right)\zeta \left( 7 \right) + 5\zeta \left( 2 \right)\zeta \left( 3 \right)\zeta \left( 5 \right)\nonumber\\
& + \frac{3}{2}{\zeta ^2}\left( 3 \right)\zeta \left( 4 \right) - \frac{5}{2}\zeta \left( 2 \right){S_{2,6}} + \frac{{211}}{{32}}{S_{2,8}},
\end{align}
We prove identities (\ref{equ:1})-(\ref{equ:8}) in Section 4, prove identities (\ref{equ:9})-(\ref{equ:11}) in Section 5 and prove identities (\ref{equ:12})-(\ref{equ:17}) in Section 6.
\section{Main Lemmas and Theorems}
In this section, we give the main lemmas and theorems of this paper, and proofs. These lemmas and theorems will be useful in the development of the main results.
\begin{lem} For integers $a,b\geq 0$, we have
\begin{align}\label{equ:3.1}
{\zeta ^ \star }\left( {{{\left\{ 2 \right\}}_b},4,{{\left\{ 2 \right\}}_a}} \right) =& 2\bar \zeta \left( { {2a + 2b + 4}} \right) - 8{\zeta ^ \star }\left( {\overline {2b + 2},1,2a + 1} \right)\nonumber\\
& + 4{\zeta ^ \star }\left( {\overline {2b + 2},2a + 2} \right) + 4{\zeta ^ \star }\left( {\overline {2b + 3},2a + 1} \right),
\end{align}
where $\bar \zeta \left(  \cdot  \right)$ denotes the alternating Riemann zeta function defined by
\[ \bar \zeta \left( s \right) :=-\zeta(\bar s)= \sum\limits_{n = 1}^\infty  {\frac{{{{\left( { - 1} \right)}^{n - 1}}}}{{{n^s}}}} ,\;{\mathop{\Re}\nolimits} \left( s \right) \ge 1,\]
and is sometimes called the Dirichlet eta function and often designated by $\eta(s)$. It is known that $\bar \zeta(s)$  is an analytic function for $\Re(s)>0$.
\end{lem}
\pf Letting $c=4$ and $n\rightarrow \infty$ in Theorem 2.3 of reference \cite{KP2013}, by a simple calculation, we may deduce the result.\hfill$\square$
\begin{lem}(\cite{Xyz2016})\label{lem:3.2}
For positive integer $ m$ , we have
\begin{align}\label{equ:3.2}
\sum\limits_{n = 1}^\infty  {\frac{{H_n^2}}{{{n^{2m}}}}{{\left( { - 1} \right)}^{n - 1}}}
 &=-{\zeta^\star}\left( {\overline{2m},2} \right) - \left( {2m + 1} \right){\zeta^\star}\left( {\overline{2m+1},1} \right)  - {\zeta^\star}\left( {{2m+1},1} \right)
\nonumber \\
&\quad - \zeta \left( 2 \right)\bar \zeta \left( {2m} \right) - 2\sum\limits_{l = 2}^m {{\zeta^\star}\left( {{2l-1},1} \right)\bar \zeta \left( {2m + 2 - 2l} \right)} ,
\end{align}
\end{lem}
\begin{lem}(\cite{TY2012}) For any non-negative integers $p,m,n$ and real $a,b,c>0$, the following identities hold:
\begin{align}\label{equ:3.3}
\zeta _p^ \star \left( {{{\left\{ c \right\}}_m},b,{{\left\{ a \right\}}_n}} \right) = \sum\limits_{k = 0}^m {\sum\limits_{l = 0}^n {{{\left( { - 1} \right)}^{k + l}}{\zeta _p}\left( {{{\left\{ a \right\}}_l},b,{{\left\{ c \right\}}_k}} \right)\zeta _p^ \star \left( {{{\left\{ c \right\}}_{m - k}}} \right)\zeta _p^ \star \left( {{{\left\{ a \right\}}_{n - l}}} \right)} } ,
\end{align}
\begin{align}\label{equ:3.4}
&{\zeta _p}\left( {{{\left\{ a \right\}}_m},b,{{\left\{ c \right\}}_n}} \right) = \sum\limits_{k = 0}^m {\sum\limits_{l = 0}^n {{{\left( { - 1} \right)}^{k + l}}\zeta _p^ \star \left( {{{\left\{ c \right\}}_l},b,{{\left\{ a \right\}}_k}} \right){\zeta _p}\left( {{{\left\{ c \right\}}_{n - l}}} \right)\zeta _p^ \star \left( {{{\left\{ a \right\}}_{m - k}}} \right)} } .
\end{align}
\end{lem}
\begin{lem}(\cite{X2016}) For positive integers $s$ and $t$, then the following relations hold:
\begin{align}\label{equ:3.5}
{\rm Li}{_s}\left( x \right){\rm Li}{_t}\left( x \right) =& \sum\limits_{j = 1}^s {A_j^{\left( {s,t} \right)}} \sum\limits_{n = 1}^\infty  {\frac{{{H^{(j)} _n}}}{{{n^{s + t - j}}}}} {x^n} + \sum\limits_{j = 1}^t {B_j^{\left( {s,t} \right)}} \sum\limits_{n = 1}^\infty  {\frac{{{H^{(j)} _n}}}{{{n^{s + t - j}}}}} {x^n}\nonumber\\& - \left( {\sum\limits_{j = 1}^s {A_j^{\left( {s,t} \right)}}  + \sum\limits_{j = 1}^t {B_j^{\left( {s,t} \right)}} } \right){\rm Li}{_{s + t}}\left( x \right),
\end{align}
where the coefficients $A_j^{\left( {s,t} \right)}$ and $B_j^{\left( {s,t} \right)}$ are defined as
\[A_j^{\left( {s,t} \right)} = \left( {\begin{array}{*{20}{c}}
   {s + t - j - 1}  \\
   {s - j}  \\
\end{array}} \right),B_j^{\left( {s,t} \right)} = \left( {\begin{array}{*{20}{c}}
   {s + t - j - 1}  \\
   {t - j}  \\
\end{array}} \right).\]
\end{lem}
\begin{lem}(\cite{Xyz2016}) For positive integers $n,m$ and real $x\in [-1,1)$, we have
\begin{align}\label{equ:3.6}
\int\limits_0^x {{t^{n - 1}}{\rm Li}{_m}\left( t \right)dt}  =& \sum\limits_{i = 1}^{m - 1} {{{\left( { - 1} \right)}^{i - 1}}\frac{{{x^n}}}{{{n^i}}}{\rm Li}{_{m + 1 - i}}\left( x \right)} - \frac{{{{\left( { - 1} \right)}^m}}}{{{n^m}}}\left( {\sum\limits_{k = 1}^n {\frac{{{x^k}}}{k}} } \right)\nonumber\\
& + \frac{{{{\left( { - 1} \right)}^m}}}{{{n^m}}}\ln \left( {1 - x} \right)\left( {{x^n} - 1} \right) ,
\end{align}
\end{lem}
Here ${\rm Li}_p(x)$ denotes the polylogarithm function defined for $x\in [-1,1)$ by
\[{\rm Li}{_p}\left( x \right) = \sum\limits_{n = 1}^\infty  {\frac{{{x^n}}}{{{n^p}}}}, \Re (p)>1, \]
with ${\rm Li_1}=-\ln(1-x),\ x\in [-1,1).$
\begin{lem}(\cite{Xu2016,X2016}) \label{lem:3.6}
For any real $l>1$, we have
\begin{align*}
&\zeta (3l){\zeta ^2}(l) = 2S_{l(3l),l} + S_ {{l^2},3l}  - S_{3l,2l} - 2S_{l,4l} + \zeta(5l),\\
&\zeta^5(l)=5S_{l^4,l}-10S_{l^3;2l}+10S_{l^2,3l}-5S_{l,4l}+\zeta(5l),\\
&\zeta (2l){\zeta ^3}(l) = {S_{{l^3},2l}} + 3{S_{{l^2}\left( {2l} \right),l}} - 3{S_{{l^2},3l}} - 3{S_{l\left( {2l} \right),2l}} + 3{S_{l,4l}} + {S_{2l,3l}} - \zeta (5l),\\
&{\zeta ^2}(2l){\zeta ^2}(3l) = 2{S_{{{\left( {2l} \right)}^2}\left( {3l} \right),3l}} + 2{S_{\left( {2l} \right){{\left( {3l} \right)}^2},2l}} - {S_{{{\left( {2l} \right)}^2},6l}} - 4{S_{\left( {2l} \right)\left( {3l} \right),5l}} - {S_{{{\left( {3l} \right)}^2},4l}}\\
&\quad\quad \quad \quad \quad \quad  + 2{S_{2l,8l}} + 2{S_{3l,7l}} - \zeta \left( {10l} \right)\quad \left( {l \ge 1} \right).
\end{align*}
\end{lem}
\pf These results immediately follows from Theorem 2.1 of \cite{Xu2016} or Theorem 2.8 of \cite{X2016}.\hfill$\square$

\begin{thm}\label{thm:1} For reals $l_1,l_2,l_3>1$ and $l_4>1$, then the following identities hold:
\begin{align}\label{equ:3.7}
\zeta \left( {{l_1}} \right)\zeta \left( {{l_2}} \right)\zeta \left( {{l_3}} \right) =& {S_{{l_1}{l_2},{l_3}}} + {S_{{l_1}{l_3},{l_2}}} + {S_{{l_2}{l_3},{l_1}}} - {S_{{l_1},{l_2} + {l_3}}} - {S_{{l_2},{l_1} + {l_3}}} - {S_{{l_3},{l_1} + {l_2}}} + \zeta \left( {{l_1} + {l_2} + {l_3}} \right),
\end{align}
\begin{align}\label{equ:3.8}
\zeta \left( {{l_1}} \right)\zeta \left( {{l_2}} \right)\zeta \left( {{l_3}} \right)\zeta \left( {{l_4}} \right) = & {S_{{l_1}{l_2}{l_3},{l_4}}} + {S_{{l_1}{l_2}{l_4},{l_3}}} + {S_{{l_1}{l_3}{l_4},{l_2}}} + {S_{{l_2}{l_3}{l_4},{l_1}}} - {S_{{l_1}{l_2},{l_3} + {l_4}}} - {S_{{l_1}{l_3},{l_2} + {l_4}}}\nonumber\\
& - {S_{{l_1}{l_4},{l_2} + {l_3}}} - {S_{{l_2}{l_3},{l_1} + {l_4}}} - {S_{{l_2}{l_4},{l_1} + {l_3}}} - {S_{{l_3}{l_4},{l_1} + {l_2}}} + {S_{{l_1},{l_2} + {l_3} + {l_4}}}\nonumber\\
& + {S_{{l_2},{l_1} + {l_3} + {l_4}}} + {S_{{l_3},{l_1} + {l_2} + {l_4}}} + {S_{{l_4},{l_1} + {l_2} + {l_3}}} - \zeta \left( {{l_1} + {l_2} + {l_3} + {l_4}} \right).
\end{align}
\end{thm}
\pf Similarly as in the proof of Theorem 2.8 of \cite{X2016}, considering the following generating functions
\begin{align*}
&{F_1}\left( x \right) := \sum\limits_{n = 1}^\infty  {\left( {H_n^{\left( {{l_1}} \right)}H_n^{\left( {{l_2}} \right)}H_n^{\left( {{l_3}} \right)} - H_n^{\left( {{l_1} + {l_2} + {l_3}} \right)}} \right){x^n}},\\
&{F_2}\left( x \right) := \sum\limits_{n = 1}^\infty  {\left( {H_n^{\left( {{l_1}} \right)}H_n^{\left( {{l_2}} \right)}H_n^{\left( {{l_3}} \right)}H_n^{\left( {{l_4}} \right)} - H_n^{\left( {{l_1} + {l_2} + {l_3} + {l_4}} \right)}} \right){x^n}} ,
\end{align*}
where $x\in(-1,1)$. Then using the definition of the harmonic numbers, we deduce that
\begin{align*}
\left( {1 - x} \right){F_1}\left( x \right) =& \sum\limits_{n = 1}^\infty  {\left\{ {\frac{{H_n^{\left( {{l_1}} \right)}H_n^{\left( {{l_2}} \right)}}}{{{n^{{l_3}}}}} + \frac{{H_n^{\left( {{l_2}} \right)}H_n^{\left( {{l_3}} \right)}}}{{{n^{{l_1}}}}} + \frac{{H_n^{\left( {{l_1}} \right)}H_n^{\left( {{l_3}} \right)}}}{{{n^{{l_2}}}}} - \frac{{H_n^{\left( {{l_1}} \right)}}}{{{n^{{l_2} + {l_3}}}}} - \frac{{H_n^{\left( {{l_2}} \right)}}}{{{n^{{l_1} + {l_3}}}}} - \frac{{H_n^{\left( {{l_3}} \right)}}}{{{n^{{l_1} + {l_2}}}}}} \right\}{x^n}} ,\\
\left( {1 - x} \right){F_2}\left( x \right) =& \sum\limits_{n = 1}^\infty  {\left\{ {\frac{{H_n^{\left( {{l_1}} \right)}H_n^{\left( {{l_2}} \right)}H_n^{\left( {{l_3}} \right)}}}{{{n^{{l_4}}}}} + \frac{{H_n^{\left( {{l_1}} \right)}H_n^{\left( {{l_2}} \right)}H_n^{\left( {{l_4}} \right)}}}{{{n^{{l_3}}}}} + \frac{{H_n^{\left( {{l_1}} \right)}H_n^{\left( {{l_3}} \right)}H_n^{\left( {{l_4}} \right)}}}{{{n^{{l_2}}}}} + \frac{{H_n^{\left( {{l_2}} \right)}H_n^{\left( {{l_3}} \right)}H_n^{\left( {{l_4}} \right)}}}{{{n^{{l_1}}}}}} \right\}{x^n}} \\
&- \sum\limits_{n = 1}^\infty  {\left\{ {\frac{{H_n^{\left( {{l_1}} \right)}H_n^{\left( {{l_2}} \right)}}}{{{n^{{l_3} + {l_4}}}}} + \frac{{H_n^{\left( {{l_1}} \right)}H_n^{\left( {{l_3}} \right)}}}{{{n^{{l_2} + {l_4}}}}} + \frac{{H_n^{\left( {{l_1}} \right)}H_n^{\left( {{l_4}} \right)}}}{{{n^{{l_2} + {l_3}}}}} + \frac{{H_n^{\left( {{l_2}} \right)}H_n^{\left( {{l_3}} \right)}}}{{{n^{{l_1} + {l_4}}}}} + \frac{{H_n^{\left( {{l_2}} \right)}H_n^{\left( {{l_4}} \right)}}}{{{n^{{l_1} + {l_3}}}}} + \frac{{H_n^{\left( {{l_3}} \right)}H_n^{\left( {{l_4}} \right)}}}{{{n^{{l_1} + {l_2}}}}}} \right\}{x^n}} \\
& + \sum\limits_{n = 1}^\infty  {\left\{ {\frac{{H_n^{\left( {{l_1}} \right)}}}{{{n^{{l_2} + {l_3} + {l_4}}}}} + \frac{{H_n^{\left( {{l_2}} \right)}}}{{{n^{{l_1} + {l_3} + {l_4}}}}} + \frac{{H_n^{\left( {{l_3}} \right)}}}{{{n^{{l_1} + {l_2} + {l_4}}}}} + \frac{{H_n^{\left( {{l_4}} \right)}}}{{{n^{{l_1} + {l_2} + {l_3}}}}}} \right\}{x^n}}  - 2\sum\limits_{n = 1}^\infty  {\frac{{{x^n}}}{{{n^{{l_1} + {l_2} + {l_3} + {l_4}}}}}} .
\end{align*}
Multiplying above equations by $(1-x)^{-1}$ and using the Cauchy product of power series, then comparing the coefficients of $x^n$ and letting $n$ tend to infinity, we may easily deduce the desired results. This completes the proof of Theorem 3.7. \hfill$\square$
\begin{re} Based on the results of Theorem \ref{thm:1}, we conjectured that for any reals $k_j>1$, the following relation holds
\begin{align}\label{equ:3.9}
\prod\limits_{j = 1}^n {\zeta \left( {{k_j}} \right)}  = \sum\limits_{j = 0}^{n - 1} {\sum\limits_{1 \le {i_1} <  \cdots  < {i_j} \le n} {{{\left( { - 1} \right)}^{n - 1 - j}}{S_{{k_{{i_1}}}{k_{{i_2}}} \cdots {k_{{i_j}}},({k_1} +  \cdots  + {k_n}) - ({k_{{i_1}}} +  \cdots  + {k_{{i_j}}})}}} } .
\end{align}
We have been proved that when $n\leq 5$, formula (\ref{equ:3.9}) is true.
It is clear that (\ref{equ:3.7}) and (\ref{equ:3.8}) are immediate corollaries of (\ref{equ:3.9}).
\end{re}
\begin{thm} For positive integers $m,p$ and real $a,b,x$ with the conditions $\left| a \right| = \left| b \right| > 1$ and $\left| x \right| < {\left| a \right|^{ - 1}}$ or ${\left| b \right|^{ - 1}}$, we have
\begin{align}\label{equ:3.10}
&{\left( { - 1} \right)^p}\sum\limits_{n = 1}^\infty  {\frac{{H_n^{\left( 1 \right)}\left( {ax} \right)H_n^{\left( m \right)}\left( b \right)}}{{{a^n}{n^p}}}}  - {\left( { - 1} \right)^m}\sum\limits_{n = 1}^\infty  {\frac{{H_n^{\left( 1 \right)}\left( {bx} \right)H_n^{\left( p \right)}\left( a \right)}}{{{b^n}{n^m}}}}  \nonumber\\
 = &\sum\limits_{i = 1}^{p - 1} {{{\left( { - 1} \right)}^{i - 1}}{\rm{L}}{{\rm{i}}_{p + 1 - i}}\left( {ax} \right)\sum\limits_{n = 1}^\infty  {\frac{{H_n^{\left( m \right)}\left( b \right)}}{{{n^i}}}{x^n}} }  - \sum\limits_{i = 1}^{m - 1} {{{\left( { - 1} \right)}^{i - 1}}{\rm{L}}{{\rm{i}}_{m + 1 - i}}\left( {bx} \right)\sum\limits_{n = 1}^\infty  {\frac{{H_n^{\left( p \right)}\left( a \right)}}{{{n^i}}}{x^n}} } \nonumber\\
 & + {\left( { - 1} \right)^p}\ln \left( {1 - ax} \right)\sum\limits_{n = 1}^\infty  {\frac{{H_n^{\left( m \right)}\left( b \right)}}{{{n^p}}}\left( {{x^n} - {a^{ - n}}} \right)}  - {\left( { - 1} \right)^m}\ln \left( {1 - bx} \right)\sum\limits_{n = 1}^\infty  {\frac{{H_n^{\left( p \right)}\left( b \right)}}{{{n^m}}}\left( {{x^n} - {b^{ - n}}} \right)},
\end{align}
where the notation ${H_n^{\left( p \right)}\left( a \right)}$ is defined by the finite sums
\[H_n^{\left( p \right)}\left( a \right) = \sum\limits_{j = 1}^n {\frac{{{a^j}}}{{{j^p}}}} .\]
If $a=1$, then $H_n^{\left( p \right)}\left( 1 \right) = H_n^{\left( p \right)}$ reduce to the harmonic numbers.
\end{thm}
\pf Applying the definition of polylogarithm function ${\rm Li}_p(x)$ and using the Cauchy product of power series, we can find that
\[\frac{{{\rm{L}}{{\rm{i}}_p}\left( {at} \right)}}{{1 - t}} = \sum\limits_{n = 1}^\infty  {H_n^{\left( p \right)}\left( a \right){t^n}} ,\;t \in \left( { - 1,1} \right).\]
Then using above identity it is easily see that
\begin{align*}
\int\limits_0^x {\frac{{{\rm{L}}{{\rm{i}}_p}\left( {at} \right){\rm{L}}{{\rm{i}}_m}\left( {bt} \right)}}{{t\left( {1 - t} \right)}}} dt =& \sum\limits_{n = 1}^\infty  {H_n^{\left( p \right)}\left( a \right)\int\limits_0^x {{t^{n - 1}}{\rm{L}}{{\rm{i}}_m}\left( {bt} \right)} dt} \\
 =& \sum\limits_{n = 1}^\infty  {H_n^{\left( m \right)}\left( b \right)\int\limits_0^x {{t^{n - 1}}{\rm{L}}{{\rm{i}}_p}\left( {at} \right)} dt} .
\end{align*}
Then with the help of formula (\ref{equ:3.6}) and the elementary transform
\[\int\limits_0^x {{t^{n - 1}}{\rm{L}}{{\rm{i}}_m}\left( {bt} \right)dt}  = \frac{1}{{{b^n}}}\int\limits_0^{bx} {{t^{n - 1}}{\rm{L}}{{\rm{i}}_m}\left( t \right)dt} \]
by a direct calculation, we may deduce the result. \hfill$\square$

\begin{thm}
For positive integers $m,p$ and real $a,b,x$ with the conditions $\left| a \right| = \left| b \right| > 1$ and $\left| x \right| < {\left| a \right|^{ - 1}}$ or ${\left| b \right|^{ - 1}}$, then the following identity holds
\begin{align}\label{equ:3.11}
&{\left( { - 1} \right)^p}\sum\limits_{n = 1}^\infty  {\frac{{H_n^{\left( 2 \right)}\left( {ax} \right)H_n^{\left( m \right)}\left( b \right)}}{{{a^n}{n^p}}}}  - {\left( { - 1} \right)^m}\sum\limits_{n = 1}^\infty  {\frac{{H_n^{\left( 2 \right)}\left( {bx} \right)H_n^{\left( p \right)}\left( a \right)}}{{{b^n}{n^m}}}} \nonumber\\
 =& \sum\limits_{i = 1}^{p - 1} {\sum\limits_{j = 1}^{p - i} {{{\left( { - 1} \right)}^{i + j}}} {\rm{L}}{{\rm{i}}_{p + 2 - i - j}}\left( {ax} \right){\rm{L}}{{\rm{i}}_{i + j,m}}\left( {x,b} \right)} - \sum\limits_{i = 1}^{m - 1} {\sum\limits_{j = 1}^{m - i} {{{\left( { - 1} \right)}^{i + j}}} {\rm{L}}{{\rm{i}}_{m + 2 - i - j}}\left( {bx} \right){\rm{L}}{{\rm{i}}_{i + j,p}}\left( {x,a} \right)} \nonumber\\
 & + {\left( { - 1} \right)^p}p\ln \left( {1 - ax} \right)\left( {{\rm{L}}{{\rm{i}}_{p + 1,m}}\left( {x,b} \right) - {\rm{L}}{{\rm{i}}_{p + 1,m}}\left( {{a^{ - 1}},b} \right)} \right)\nonumber\\
 & - {\left( { - 1} \right)^m}m\ln \left( {1 - bx} \right)\left( {{\rm{L}}{{\rm{i}}_{m + 1,p}}\left( {x,a} \right) - {\rm{L}}{{\rm{i}}_{m + 1,p}}\left( {{b^{ - 1}},a} \right)} \right)\nonumber\\
 & - {\left( { - 1} \right)^p}p\sum\limits_{n = 1}^\infty  {\frac{{H_n^{\left( 1 \right)}\left( {ax} \right)H_n^{\left( m \right)}\left( b \right)}}{{{a^n}{n^{p + 1}}}}}  + {\left( { - 1} \right)^m}m\sum\limits_{n = 1}^\infty  {\frac{{H_n^{\left( 1 \right)}\left( {bx} \right)H_n^{\left( p \right)}\left( a \right)}}{{{b^n}{n^{m + 1}}}}}\nonumber \\
 & + {\left( { - 1} \right)^p}{\rm{L}}{{\rm{i}}_2}\left( {ax} \right){\rm{L}}{{\rm{i}}_{p,m}}\left( {{a^{ - 1}},b} \right) - {\left( { - 1} \right)^m}{\rm{L}}{{\rm{i}}_2}\left( {bx} \right){\rm{L}}{{\rm{i}}_{m,p}}\left( {{b^{ - 1}},a} \right),
\end{align}
where ${\rm{L}}{{\rm{i}}_{p,m}}\left( {x,y} \right)$ is defined by the double sum
\[{\rm{L}}{{\rm{i}}_{s,t}}\left( {x,y} \right): = \sum\limits_{n = 1}^\infty  {\frac{{{x^n}}}{{{n^s}}}\sum\limits_{k = 1}^n {\frac{{{y^k}}}{{{k^t}}}} } ,\;\left| {x} \right| \leq 1,\left| {y} \right| > 0,\left| {xy} \right| < 1,s,t>0\]
Of course, if $s>1$, then we can allow $\left| {xy} \right|=1$.
\end{thm}
\pf Replacing $x$ by $t$ in (\ref{equ:3.10}), then dividing it by $t$ and integrating over the interval $(0,x)$ with the help of formulae (\ref{equ:3.6}) and
\[\int\limits_0^x {\frac{{\ln \left( {1 - at} \right)}}{t}dt}  =  - {\rm{L}}{{\rm{i}}_2}\left( {ax} \right),\]
by a simple calculation, we have the result.\hfill$\square$

\begin{thm} \label{thm:3.10}For any non-negative integers $b$ and $n$, then
\begin{align}\label{equ:3.12}
&\zeta \left( {4,{{\left\{ 2 \right\}}_n}} \right) = \sum\limits_{l = 0}^n {{{\left( { - 1} \right)}^l}{\zeta ^ \star }\left( {{{\left\{ 2 \right\}}_l},4} \right)\zeta \left( {{{\left\{ 2 \right\}}_{n - l}}} \right)} ,
\end{align}
\begin{align}\label{equ:3.13}
&{\zeta ^ \star }\left( {{{\left\{ 2 \right\}}_b},4} \right) = 2\bar \zeta \left( {2b + 4} \right) + 4\sum\limits_{n = 1}^\infty  {\frac{{H_n^2}}{{{n^{2b + 2}}}}{{\left( { - 1} \right)}^{n - 1}}}  - 4\sum\limits_{n = 1}^\infty  {\frac{{{H_n}}}{{{n^{2b + 3}}}}{{\left( { - 1} \right)}^{n - 1}}} ,
\end{align}
where \[\zeta \left( {{{\left\{ 2 \right\}}_b}} \right) = \frac{{{\pi ^{2b}}}}{{\left( {2b + 1} \right)!}}.\]
\end{thm}
\pf Letting $m=0,b=4,c=2,a=0$ in (\ref{equ:3.4}) and $p\rightarrow \infty$, we obtain (\ref{equ:3.12}).
Taking $a=0$ in (\ref{equ:3.1}) and using the well-known fact
\[{\zeta _{n - 1}}\left( {1,1} \right) = \frac{{H_n^2 - H_n^{\left( 2 \right)}}}{2} - \frac{{{H_n}}}{n} + \frac{1}{{{n^2}}}\]
the result is formula (\ref{equ:3.13}). \hfill$\square$\\
From Lemma \ref{lem:3.2} and Theorem \ref{thm:3.10}, we know that for integer $b\geq0$, then the MZVs $\zeta \left( {4,{{\left\{ 2 \right\}}_b}} \right)$ and MZSVs ${\zeta ^ \star }\left( {{{\left\{ 2 \right\}}_b},4} \right)$ are reducible to linear sums. For example, from (\ref{equ:3.2}), (\ref{equ:3.12}) and (\ref{equ:3.13}), we get
\begin{align}\label{equ:3.14}
\zeta \left( {4,{{\left\{ 2 \right\}}_3}} \right) =& \frac{{20719}}{{256}}\zeta \left( {10} \right) - 2{\zeta ^2}\left( 5 \right) - 4\zeta \left( 3 \right)\zeta \left( 7 \right) - 24\zeta \left( 2 \right)\zeta \left( 3 \right)\zeta \left( 5 \right)\nonumber\\
& + \frac{9}{4}{\zeta ^2}\left( 3 \right)\zeta \left( 4 \right) + \frac{9}{2}\zeta \left( 2 \right){S_{2,6}} - 4{{\bar S}_{2,8}} - 32{{\bar S}_{1,9}},
\end{align}

\begin{align}\label{equ:3.15}
{\zeta ^ \star }\left( {{{\left\{ 2 \right\}}_3},4} \right) =&  - \frac{{69803}}{{1280}}\zeta \left( {10} \right) + 2{\zeta ^2}\left( 5 \right) + 4\zeta \left( 3 \right)\zeta \left( 7 \right) + 4\zeta \left( 2 \right)\zeta \left( 3 \right)\zeta \left( 5 \right)\nonumber\\
& + \frac{7}{2}{\zeta ^2}\left( 3 \right)\zeta \left( 4 \right) + 4{{\bar S}_{2,8}} + 32{{\bar S}_{1,9}},
\end{align}

\begin{align}\label{equ:3.16}
{\zeta ^ \star }\left( {2,2,4} \right) =  - \frac{{6787}}
{{192}}\zeta \left( 8 \right) + 4\zeta \left( 3 \right)\zeta \left( 5 \right) + 2\zeta \left( 2 \right){\zeta ^2}\left( 3 \right) + 4{{\bar S}_{2,6}} + 24{{\bar S}_{1,7}},
\end{align}

\begin{align}\label{equ:3.17}
\zeta \left( {4,2,2} \right) =  - \frac{{22817}}
{{576}}\zeta \left( 8 \right) + 4\zeta \left( 3 \right)\zeta \left( 5 \right) + 3\zeta \left( 2 \right){\zeta ^2}\left( 3 \right) + 4{{\bar S}_{2,6}} + 24{{\bar S}_{1,7}},
\end{align}

\begin{align}\label{equ:3.18}
\sum\limits_{n = 1}^\infty  {\frac{{H_n^2}}{{{n^8}}}{{\left( { - 1} \right)}^{n - 1}}}  &=  - \frac{{36179}}{{2560}}\zeta \left( {10} \right) + \frac{1}{2}{\zeta ^2}\left( 5 \right) + \zeta \left( 3 \right)\zeta \left( 7 \right) + \zeta \left( 2 \right)\zeta \left( 3 \right)\zeta \left( 5 \right)\nonumber\\
&\quad + \frac{7}{8}{\zeta ^2}\left( 3 \right)\zeta \left( 4 \right) + {{\bar S}_{2,8}} + 9{{\bar S}_{1,9}}.
\end{align}

On the other hand, from \cite{FS1998}, we deduce that
\begin{align}\label{equ:3.19}
\zeta \left( {4,{{\left\{ 2 \right\}}_2}} \right) =& \frac{{947}}{{72}}\zeta \left( 8 \right) + 3\zeta \left( 2 \right){\zeta ^2}\left( 3 \right) - 20\zeta \left( 3 \right)\zeta \left( 5 \right) + \frac{9}{2}{S_{2,6}},
\end{align}
Hence, from formulas (\ref{equ:3.17}) and (\ref{equ:3.19}), we obtain the following relation 
\[{{\bar S}_{2,6}} + 6{{\bar S}_{1,7}} = \frac{{3377}}
{{256}}\zeta \left( 8 \right) - 6\zeta \left( 3 \right)\zeta \left( 5 \right) + \frac{9}
{8}{S_{2,6}}.\]

\begin{cor} For positive integers $m>1$ and $p>1$, we have
\begin{align}\label{equ:3.20}
&{\left( { - 1} \right)^p}{S_{2m,p}} - {\left( { - 1} \right)^m}{S_{2p,m}}\nonumber \\=& \sum\limits_{i = 1}^{p - 1} {\sum\limits_{j = 1}^{p - i} {{{\left( { - 1} \right)}^{i + j}}\zeta \left( {p + 2 - i - j} \right)} {S_{m,i + j}}}\nonumber \\
& - \sum\limits_{i = 1}^{m - 1} {\sum\limits_{j = 1}^{m - i} {{{\left( { - 1} \right)}^{i + j}}\zeta \left( {m + 2 - i - j} \right)} {S_{p,i + j}}} \nonumber\\
& - {\left( { - 1} \right)^p}p{S_{1m,p + 1}} + {\left( { - 1} \right)^m}m{S_{1p,m + 1}}\nonumber\\
& + {\left( { - 1} \right)^p}\zeta \left( 2 \right){S_{m,p}} - {\left( { - 1} \right)^m}\zeta \left( 2 \right){S_{p,m}}.
\end{align}
\end{cor}
\pf Letting $a,b,x\rightarrow 1$ in (\ref{equ:3.11}) we obtain the result. \hfill$\square$
\section{Proofs of identities (\ref{equ:1})-(\ref{equ:8})}
In this section we prove that equations (\ref{equ:1})-(\ref{equ:8}) hold by using the methods of integrals of polylogarithms and above lemmas and theorems.

 We begin with some well-known identities. In \cite{X2017}, we shown that for weight $1+p+q=10$, all quadratic sums $S_{1p_,q}$ are reducible to $S_{2,6}$ and $S_{2,8}$. The explicit formulas are as follows
 \begin{align}\label{4.1}
S_{12,7}
=&  - \frac{{1331}}{{80}}\zeta( {10}) + \frac{{43}}{4}{\zeta ^2}( 5) + \frac{{41}}{2}\zeta( 3)\zeta( 7)
- 7\zeta( 2)\zeta( 3)\zeta( 5)\nonumber \\&- 2{\zeta ^2}( 3)\zeta( 4)
- \frac{5}{4}{S_{2,8}},
\end{align}
 \begin{align}\label{4.2}
S_{13,6} =&  - \frac{{247}}
{{40}}\zeta( {10}) - \frac{5}
{4}{\zeta ^2}( 5) - \frac{{15}}
{2}\zeta( 3)\zeta( 7)\nonumber\\
& + 12\zeta( 2)\zeta( 3)\zeta( 5) - \frac{{21}}
{4}{S_{2,8}} - \zeta( 2){S_{2,6}},
\end{align}
\begin{align}\label{4.3}
S_{2^2,6}  =& \frac{{2697}}
{{40}}\zeta( {10}) - 41{\zeta ^2}( 5) - 63\zeta( 3)\zeta( 7) + 16\zeta( 2)\zeta( 3)\zeta( 5)\nonumber\\
& + 4{\zeta ^2}( 3)\zeta( 4) + \frac{{23}}
{2}{S_{2,8}}{\text{ + }}2\zeta( 2){S_{2,6}},
\end{align}
\begin{align}\label{4.4}
S_{14,5} = &\frac{{6033}}
{{160}}\zeta( {10}) - 14{\zeta ^2}( 5) - 4\zeta( 3)\zeta( 7) - 15\zeta( 2 )\zeta( 3)\zeta( 5)\nonumber\\
& - \frac{1}
{2}{\zeta ^2}( 3)\zeta( 4) + \frac{{21}}
{2}{S_{2,8}} + \frac{5}
{2}\zeta( 2){S_{2,6}},
\end{align}
\begin{align}\label{4.5}
S_{15,4} = & - \frac{{6569}}
{{240}}\zeta( {10}) + 16{\zeta ^2}( 5) + 10\zeta( 3)\zeta( 7) + 4\zeta( 2)\zeta( 3)\zeta( 5)\nonumber\\
& + {\zeta ^2}( 3)\zeta( 4) - \frac{{21}}
{2}{S_{2,8}},
\end{align}
\begin{align}\label{4.6}
S_{16,3}  =& \frac{{1043}}
{{160}}\zeta( {10}) - \frac{{17}}
{4}{\zeta ^2}( 5) - \frac{{15}}
{2}\zeta( 3)\zeta( 7) + 4\zeta( 2)\zeta( 3)\zeta( 5)\nonumber\\
& - \frac{1}
{2}{\zeta ^2}( 3)\zeta( 4) - \frac{5}
{2}\zeta( 2){S_{2,6}} + \frac{{21}}
{4}{S_{2,8}},
\end{align}
\begin{align}\label{4.7}
S_{17,2}  =& \frac{{242}}
{{15}}\zeta( {10}) - \frac{{25}}
{4}{\zeta ^2}( 5) - \frac{{19}}
{2}\zeta( 3)\zeta( 7)+{\zeta ^2}( 3)\zeta( 4)\nonumber\\
&+\zeta( 2){S_{2,6}} + \frac{5}
{4}{S_{2,8}}.
\end{align}
From lemma \ref{lem:3.6} and theorem \ref{thm:1}, the following identities are easily derived
\begin{align}\label{4.8}
&{\zeta ^2}\left( 2 \right)\zeta \left( 6 \right) = {S_{{2^2},6}} + 2{S_{26,2}} - 2{S_{2,8}} - {S_{6,4}} + \zeta \left( {10} \right),
\end{align}
\begin{align}\label{4.9}
&{\zeta ^2}\left( 4 \right)\zeta \left( 2 \right) = {S_{{4^2},2}} + 2{S_{24,4}} - 2{S_{4,6}} - {S_{2,8}} + \zeta \left( {10} \right),
\end{align}
\begin{align}\label{4.10}
&{\zeta ^2}\left( 3 \right)\zeta \left( 4 \right) = {S_{{3^2},4}} + 2{S_{34,3}} - 2{S_{3,7}} - {S_{4,6}} + \zeta \left( {10} \right),
\end{align}
\begin{align}\label{4.11}
&\zeta \left( 2 \right)\zeta \left( 3 \right)\zeta \left( 5 \right) = {S_{23,5}} + {S_{25,3}} + {S_{35,2}} - {S_{5,5}} - {S_{2,8}} - {S_{3,7}} + \zeta \left( {10} \right).
\end{align}
Taking $m=3,p=5$ in (\ref{equ:3.20}) yields
\begin{align}\label{4.12}
 {S_{25,3}} - {S_{23,5}} = &\frac{{7179}}{{80}}\zeta \left( {10} \right) - \frac{{195}}{4}{\zeta ^2}\left( 5 \right) - \frac{{49}}{2}\zeta \left( 3 \right)\zeta \left( 7 \right) - 17\zeta \left( 2 \right)\zeta \left( 3 \right)\zeta \left( 5 \right)\nonumber\\
 &- 2{\zeta ^2}\left( 3 \right)\zeta \left( 4 \right) + 15\zeta \left( 2 \right){S_{2,6}} + \frac{{21}}{4}{S_{2,8}}.
\end{align}
From formula (3.1) in the reference \cite{X2017}, letting $s=t=2,r=5$, we obtain
\begin{align}\label{4.13}
&8{S_{12,7}} + 6{S_{13,6}} + 4{S_{14,5}} + 2{S_{15,4}} + \frac{3}{2}{S_{{2^2},6}} + 2{S_{23,5}} + {S_{24,4}} + \frac{1}{2}{S_{{3^2},4}}\nonumber\\
 =&  - \frac{{6301}}{{120}}\zeta \left( {10} \right) + 39{\zeta ^2}\left( 5 \right) + 88\zeta \left( 3 \right)\zeta \left( 7 \right) - 22\zeta \left( 2 \right)\zeta \left( 3 \right)\zeta \left( 5 \right)\nonumber\\
 & - 9{\zeta ^2}\left( 3 \right)\zeta \left( 4 \right) + 2\zeta \left( 2 \right){S_{2,6}} - 14{S_{2,8}}.
\end{align}
Next, we consider the following two integrals of polylogarithms
\[\int\limits_0^1 {\frac{{{\rm{Li}}_3^3\left( x \right)}}{x}} dx\quad {\rm and}\quad \int\limits_0^1 {\frac{{{\rm{Li}}_2^2\left( x \right){\rm{L}}{{\rm{i}}_4}\left( x \right)\ln x}}{x}} dx.\]
On the one hand, setting $s=t=3$ in (\ref{equ:3.5}) we deduce that
\begin{align}\label{equ:4.14}
{\rm{Li}}_3^2\left( x \right) = 2\sum\limits_{n = 1}^\infty  {\left\{ {\frac{{H_n^{\left( 3 \right)}}}{{{n^3}}} + 3\frac{{H_n^{\left( 2 \right)}}}{{{n^4}}} + 6\frac{{{H_n}}}{{{n^5}}} - 10\frac{1}{{{n^6}}}} \right\}{x^n}} .
\end{align}
Multiplying (\ref{equ:4.14}) by $\frac{{{\rm  Li}_2(x)}}{x}$ and integrating over $(0,x)$ with the help of formula (\ref{equ:3.6}). The result is
\begin{align}\label{equ:4.15}
{\rm{Li}}_3^3\left( x \right) = 6\sum\limits_{n = 1}^\infty  {\left\{ {\frac{{H_n^{\left( 3 \right)}}}{{{n^3}}} + 3\frac{{H_n^{\left( 2 \right)}}}{{{n^4}}} + 6\frac{{{H_n}}}{{{n^5}}} - 10\frac{1}{{{n^6}}}} \right\}\left\{ {\frac{{{x^n}}}{n}{\rm{L}}{{\rm{i}}_2}\left( x \right) + \frac{{{x^n} - 1}}{{{n^2}}}\ln \left( {1 - x} \right) - \frac{1}{{{n^2}}}\left( {\sum\limits_{k = 1}^n {\frac{{{x^k}}}{k}} } \right)} \right\}}.
\end{align}
Hence, by a direct calculation, we have
\begin{align}\label{4.16}
\int\limits_0^1 {\frac{{{\rm{Li}}_3^3\left( x \right)}}{x}} dx =& 6\sum\limits_{n = 1}^\infty  {\left\{ {\frac{{H_n^{\left( 3 \right)}}}{{{n^3}}} + 3\frac{{H_n^{\left( 2 \right)}}}{{{n^4}}} + 6\frac{{{H_n}}}{{{n^5}}} - 10\frac{1}{{{n^6}}}} \right\}\left\{ {2\frac{{\zeta \left( 2 \right)}}{{{n^2}}} - 2\frac{{{H_n}}}{{{n^3}}} - \frac{{H_n^{\left( 2 \right)}}}{{{n^2}}}} \right\}}\nonumber\\
=&  - \frac{{17829}}{{40}}\zeta \left( {10} \right) + 279{\zeta ^2}\left( 5 \right) + 348\zeta \left( 3 \right)\zeta \left( 7 \right) - 84\zeta \left( 2 \right)\zeta \left( 3 \right)\zeta \left( 5 \right)\nonumber\\
&- 18\zeta \left( 2 \right){S_{2,6}} - 66{S_{2,8}} - 6{S_{23,5}}.
\end{align}
On the other hand, multiplying (\ref{equ:4.14}) by $\frac{{{\rm  Li}_3(x)}}{x}$ and integrating over $(0,1)$, we obtain \begin{align}\label{4.17}
\int\limits_0^1 {\frac{{{\rm{Li}}_3^3\left( x \right)}}{x}} dx =& 2\sum\limits_{n = 1}^\infty  {\left\{ {\frac{{H_n^{\left( 3 \right)}}}{{{n^3}}} + 3\frac{{H_n^{\left( 2 \right)}}}{{{n^4}}} + 6\frac{{{H_n}}}{{{n^5}}} - 10\frac{1}{{{n^6}}}} \right\}\left\{ {\frac{{\zeta \left( 3 \right)}}{n} - \frac{{\zeta \left( 2 \right)}}{{{n^2}}} + \frac{{{H_n}}}{{{n^3}}}} \right\}}\nonumber \\
 = & - \frac{{2253}}{{80}}\zeta \left( {10} \right) + 12{\zeta ^2}\left( 5 \right) + 12\zeta \left( 3 \right)\zeta \left( 7 \right) + 6\left( 2 \right)\zeta \left( 3 \right)\zeta \left( 5 \right)\nonumber\\
 &- 3\zeta \left( 2 \right){S_{2,6}} - 6{S_{2,8}}.
\end{align}
The relations (\ref{4.16}) and (\ref{4.17}) yield the formula (\ref{equ:2}).

Similarly, letting $(s,t)=(2,2)$ and $(2,4)$ in (\ref{equ:3.5}) we can find that
\begin{align*}
&{\rm{Li}}_2^2\left( x \right) = 2\sum\limits_{n = 1}^\infty  {\left\{ {\frac{{H_n^{\left( 2 \right)}}}{{{n^2}}} + 2\frac{{{H_n}}}{{{n^3}}} - 3\frac{1}{{{n^4}}}} \right\}{x^n}} ,\\
&{\rm{L}}{{\rm{i}}_2}\left( x \right){\rm{L}}{{\rm{i}}_4}\left( x \right) = \sum\limits_{n = 1}^\infty  {\left\{ {\frac{{H_n^{\left( 4 \right)}}}{{{n^2}}} + 2\frac{{H_n^{\left( 3 \right)}}}{{{n^3}}} + 4\frac{{H_n^{\left( 2 \right)}}}{{{n^4}}} + 8\frac{{{H_n}}}{{{n^5}}} - 15\frac{1}{{{n^6}}}} \right\}{x^n}} .
\end{align*}
Dividing (\ref{equ:3.6}) by $x$ and integrating over the interval $(0,1)$, by a simple calculation, we arrive at the conclusion that
\begin{align}\label{4.18}
\int\limits_0^1 {{x^{n - 1}}\ln x{\rm{L}}{{\rm{i}}_m}\left( x \right)} dx = &\sum\limits_{i = 1}^{m - 1} {\sum\limits_{j = 1}^{m - i} {{{\left( { - 1} \right)}^{i + j - 1}}\frac{{\zeta \left( {m + 2 - i - j} \right)}}{{{n^{i + j}}}}} }\nonumber \\
&+ {\left( { - 1} \right)^m}m\frac{{{H_n}}}{{{n^{m + 1}}}} + {\left( { - 1} \right)^m}\frac{{{H^{(2)}_n} - \zeta \left( 2 \right)}}{{{n^m}}}.
\end{align}
Setting $m=2,4$ in above equation, we conclude that
\begin{align*}
&\int\limits_0^1 {{x^{n - 1}}{\rm{ln}}x{\rm{L}}{{\rm{i}}_2}\left( x \right)} dx = 2\frac{{{H_n}}}{{{n^3}}} + \frac{{H_n^{\left( 2 \right)}}}{{{n^2}}} - 2\frac{{\zeta \left( 2 \right)}}{{{n^2}}},\\
&\int\limits_0^1 {{x^{n - 1}}{\rm{ln}}x{\rm{L}}{{\rm{i}}_4}\left( x \right)} dx =  - \frac{{\zeta \left( 4 \right)}}{{{n^2}}} + 2\frac{{\zeta \left( 3 \right)}}{{{n^3}}} - 4\frac{{\zeta \left( 2 \right)}}{{{n^4}}} + 4\frac{{{H_n}}}{{{n^5}}} + \frac{{H_n^{\left( 2 \right)}}}{{{n^4}}}.
\end{align*}
Thus, we have
\begin{align}\label{4.19}
\int\limits_0^1 {\frac{{{\rm{Li}}_2^2\left( x \right){\rm{L}}{{\rm{i}}_4}\left( x \right)\ln x}}{x}} dx = & \sum\limits_{n = 1}^\infty  {\left\{ {\frac{{H_n^{\left( 4 \right)}}}{{{n^2}}} + 2\frac{{H_n^{\left( 3 \right)}}}{{{n^3}}} + 4\frac{{H_n^{\left( 2 \right)}}}{{{n^4}}} + 8\frac{{{H_n}}}{{{n^5}}} - 15\frac{1}{{{n^6}}}} \right\}\int\limits_0^1 {{x^{n - 1}}{\rm{ln}}x{\rm{L}}{{\rm{i}}_2}\left( x \right)} dx}\nonumber \\
= &\frac{{1023}}{{40}}\zeta \left( {10} \right) - {\zeta ^2}\left( 5 \right) + 20\zeta \left( 3 \right)\zeta \left( 7 \right) - 32\zeta \left( 2 \right)\zeta \left( 3 \right)\zeta \left( 5 \right)\nonumber\\
& - {\zeta ^2}\left( 3 \right)\zeta \left( 4 \right) + 6\zeta \left( 2 \right){S_{2,6}} + 7{S_{2,8}} + {S_{24,4}}\nonumber\\
 =& 2\sum\limits_{n = 1}^\infty  {\left\{ {\frac{{H_n^{\left( 2 \right)}}}{{{n^2}}} + 2\frac{{{H_n}}}{{{n^3}}} - 3\frac{1}{{{n^4}}}} \right\}\int\limits_0^1 {{x^{n - 1}}{\rm{ln}}x{\rm{L}}{{\rm{i}}_4}\left( x \right)} dx}\nonumber \\
= &\frac{{1147}}{{30}}\zeta \left( {10} \right) - 21{\zeta ^2}\left( 5 \right) - 36\zeta \left( 3 \right)\zeta \left( 7 \right) + 8\zeta \left( 2 \right)\zeta \left( 3 \right)\zeta \left( 5 \right)\nonumber\\
&- 4\zeta \left( 2 \right){S_{2,6}} + 18{S_{2,8}}.\tag{4.19}
\end{align}
Hence, combining equations (\ref{4.1})-(\ref{4.13}) and (\ref{4.19}), by a direct calculation we can obtain the closed forms of formulas (\ref{equ:1}) and (\ref{equ:3})-(\ref{equ:8}).

\section{Proofs of identities (\ref{equ:9})-(\ref{equ:11})}
From \cite{X2017,Xu2017}, we can find the relation
\[\sum\limits_{1 \le i < j < k \le n} {\frac{1}{{{i^2}{j^2}{k^2}}}}  = \frac{1}{{3!}}\left\{ {{{\left[ {H_n^{\left( 2 \right)}} \right]}^3} - 3H_n^{\left( 2 \right)}H_n^{\left( 4 \right)} + 2H_n^{\left( 6 \right)}} \right\}.\]
According to the definition of multiple zeta values, it is easily seen that
\begin{align}\label{equ:5.1}
\zeta \left( {4,2,2,2} \right) &= \frac{1}{6}\sum\limits_{n = 1}^\infty  {\frac{[{H^{(2)}  _n]^3 - 3{H^{(2)}_n}{H^{(4)}_n} + 2{H^{(6)}_n}}}{{{n^4}}}}\nonumber \\
& - \frac{1}{2}\left\{ {\sum\limits_{n = 1}^\infty  {\frac{{[H^{(2)}_n]^2}}{{{n^6}}} - 2\sum\limits_{n = 1}^\infty  {\frac{{{H^{(2)}_n}}}{{{n^8}}}}  - \sum\limits_{n = 1}^\infty  {\frac{{{H^{(4)}_n}}}{{{n^6}}}} }  + 2\zeta \left( 10 \right)} \right\}.
\end{align}
So, substituting (\ref{equ:3.14}) and (\ref{equ:5}) into (\ref{equ:5.1}) respectively, we deduce the result (\ref{equ:9}).

From lemma \ref{lem:3.6} and theorem \ref{thm:1}, we get
\begin{align}\label{equ:5.2}
{\zeta ^5}\left( 2 \right) = 5{S_{{2^4},2}} - 10{S_{{2^3},4}} + 10{S_{{2^2},6}} - 5{S_{2,8}} + \zeta \left( {10} \right),
\end{align}
\begin{align}\label{equ:5.3}
&\zeta \left( 4 \right){\zeta ^3}\left( 2 \right) = {S_{{2^3},4}} + 3{S_{{2^2}4,2}} - 3{S_{{2^2},6}} - 3{S_{24,4}} + 3{S_{2,8}} + {S_{4,6}} - \zeta \left( {10} \right).
\end{align}
Therefore, combining related equations, we can obtain the formulas (\ref{equ:10}) and (\ref{equ:11}).

In fact, by using the MZV-interface \url{http://wayback.cecm.sfu.ca/cgi-bin/EZFace/zetaform.cgi} we can find the relation
\begin{align}\label{equ:5.4}\zeta(\bar{8},2)+8\zeta(\bar{9},1)=\frac{3}{1024} \left( 1392 \zeta^2(5) + 2768 \zeta(3) \zeta(7) -4839 \zeta(10) - 376 \zeta(8,2)\right).\end{align}
Hence, by the definitions of mzvs and Euler sums, we have
\begin{align}\label{equ:5.5}
{{\bar S}_{2,8}} + 8{{\bar S}_{1,9}} = \frac{{22587}}
{{1024}}\zeta \left( {10} \right) - \frac{{261}}
{{64}}{\zeta ^2}\left( 5 \right) - \frac{{519}}
{{64}}\zeta \left( 3 \right)\zeta \left( 7 \right) + \frac{{141}}
{{128}}{S_{2,8}}.
\end{align}
Thus, using the experimental evaluation (\ref{equ:5.5}), then the identities (\ref{equ:9})-(\ref{equ:11}) can be rewritten as
\begin{align*}
 {S_{{2^3},4}} =& \frac{{29907}}{{160}}\zeta \left( {10} \right) - \frac{{737}}{8}{\zeta ^2}\left( 5 \right) - \frac{{1435}}{8}\zeta \left( 3 \right)\zeta \left( 7 \right) + 24\zeta \left( 2 \right)\zeta \left( 3 \right)\zeta \left( 5 \right) \\
  &+ \frac{{57}}{2}{\zeta ^2}\left( 3 \right)\zeta \left( 4 \right) + 3\zeta \left( 2 \right){S_{2,6}} + \frac{{505}}{{16}}{S_{2,8}}, \\
 {S_{{2^4},2}} =& \frac{{38591}}{{160}}\zeta \left( {10} \right) - \frac{{409}}{4}{\zeta ^2}\left( 5 \right) - \frac{{931}}{4}\zeta \left( 3 \right)\zeta \left( 7 \right) + 16\zeta \left( 2 \right)\zeta \left( 3 \right)\zeta \left( 5 \right) \\
  &+ 49{\zeta ^2}\left( 3 \right)\zeta \left( 4 \right) + 2\zeta \left( 2 \right){S_{2,6}} + \frac{{329}}{8}{S_{2,8}}, \\
 {S_{{2^2}4,2}} =& \frac{{2549}}{{160}}\zeta \left( {10} \right) - \frac{{229}}{8}{\zeta ^2}\left( 5 \right) - \frac{{455}}{8}\zeta \left( 3 \right)\zeta \left( 7 \right) + 48\zeta \left( 2 \right)\zeta \left( 3 \right)\zeta \left( 5 \right) \\
 & - \frac{9}{2}{\zeta ^2}\left( 3 \right)\zeta \left( 4 \right) - 9\zeta \left( 2 \right){S_{2,6}} + \frac{{157}}{{16}}{S_{2,8}}.
\end{align*}
It should be emphasized that the result (\ref{equ:5.4}) are not established in any rigorous mathematical sense. However, we can find that the error of (\ref{equ:5.4}) is better than $10^{48}$ by using EZ Face. Hence, we believe that it is correct. 

\section{Proofs of identities (\ref{equ:12})-(\ref{equ:17})}
To prove the identities (\ref{equ:12})-(\ref{equ:17}), we need the following two lemmas.
\begin{lem}(\cite{Xu2017})\label{lem:6.1}For integers $p>0$ and $m\geq 0$. Then
 \begin{align}\label{6.1}
 \sum\limits_{n = 1}^\infty  {\frac{{{H_n}s\left( {n,p} \right)}}{{n!{n^{m + 1}}}}}  = \left( {p + 1} \right)\zeta \left( {p + 2,{{\left\{ 1 \right\}}_m}} \right) + \sum\limits_{i = 1}^m {\zeta \left( {p + 1,{{\left\{ 1 \right\}}_{i - 1}},2,{{\left\{ 1 \right\}}_{m - i}}} \right)}.\end{align}
\end{lem}
\begin{lem}(\cite{Xu2017})\label{lem:6.2} For integers $m,p> 0$ and $r>1$. Then
\begin{align}\label{6.2}
\sum\limits_{n = 1}^\infty  {\frac{{s\left( {n,m} \right){H^{(r)}_n}}}{{n!{n^p}}}}  = \zeta \left( r \right)\zeta \left( {m + 1,{{\left\{ 1 \right\}}_{p - 1}}} \right) - \zeta \left( {m + 1,{{\left\{ 1 \right\}}_{p - 1}},2,{{\left\{ 1 \right\}}_{r - 2}}} \right).
\end{align}
\end{lem}
Here  ${s\left( {n,k} \right)}$ denotes the (unsigned) Stirling number of the first kind (see \cite{L1974}), and we have
\begin{align*}
& s\left( {n,1} \right) = \left( {n - 1} \right)!,\\&s\left( {n,2} \right) = \left( {n - 1} \right)!{H_{n - 1}},\\&s\left( {n,3} \right) = \frac{{\left( {n - 1} \right)!}}{2}\left[ {H_{n - 1}^2 - {H^{(2)} _{n - 1}}} \right],\\
&s\left( {n,4} \right) = \frac{{\left( {n - 1} \right)!}}{6}\left[ {H_{n - 1}^3 - 3{H_{n - 1}}{H^{(2)} _{n - 1}} + 2{H^{(3)}_{n - 1}}} \right], \\
&s\left( {n,5} \right) = \frac{{\left( {n - 1} \right)!}}{{24}}\left[ {H_{n - 1}^4 - 6{H^{(4)}_{n - 1}} - 6H_{n - 1}^2{H^{(2)}_{n - 1}} + 3(H^{(2)}_{n-1})^2+ 8H_{n - 1}^{}{H^{(3)}_{n - 1}}} \right].
\end{align*}
Noticing that Borwein, Bradley and Broadhurst \cite{BBBL1997} proved the following duality relation
\begin{align}\label{equ:6.3}
\zeta \left( {{m_1} + 2,{{\left\{ 1 \right\}}_{{n_1}}}, \ldots ,{m_p} + 2,{{\left\{ 1 \right\}}_{{n_p}}}} \right) = \zeta \left( {{n_p} + 2,{{\left\{ 1 \right\}}_{{m_p}}}, \ldots ,{n_1} + 2,{{\left\{ 1 \right\}}_{{m_1}}}} \right).
\end{align}
Hence, taking $m=3$ in Lemma \ref{lem:6.2}, we obtain
\begin{align}\label{equ:6.4}
{S_{{1^2}r,p + 1}} = {S_{2r,p + 1}} + 2{S_{1r,p + 2}} - 2{S_{r,p + 3}} + 2\zeta \left( r \right)\zeta \left( {4,{{\left\{ 1 \right\}}_{p - 1}}} \right) - 2\zeta \left( {r,p + 1,1,1} \right).
\end{align}
From \cite{Xu2017}, we can deduce the following results
\begin{align*}
&\zeta \left( {4,{{\left\{ 1 \right\}}_2}} \right) = \frac{{23}}{{16}}\zeta \left( 6 \right) - {\zeta ^2}\left( 3 \right),\\
&\zeta \left( {4,{{\left\{ 1 \right\}}_3}} \right) = 5\zeta \left( 7 \right) - \frac{5}{4}\zeta \left( 3 \right)\zeta \left( 4 \right) - 2\zeta \left( 2 \right)\zeta \left( 5 \right),\\
&\zeta \left( {4,{{\left\{ 1 \right\}}_4}} \right) = \frac{{61}}{{24}}\zeta \left( 8 \right) - 3\zeta \left( 3 \right)\zeta \left( 5 \right) + \frac{1}{2}{\zeta ^2}\left( 2 \right)\zeta \left( 3 \right).
\end{align*}
Thus, letting $(r,p)=(2,5),(3,4),(4,3),(5,2)$ and $(6,1)$ in (\ref{equ:6.4}), we have
\begin{align}\label{equ:6.5}
{S_{{1^2}2,6}} = &\frac{{3403}}{{80}}\zeta \left( {10} \right) - \frac{{39}}{2}{\zeta ^2}\left( 5 \right) - 22\zeta \left( 3 \right)\zeta \left( 7 \right) - 4\zeta \left( 2 \right)\zeta \left( 3 \right)\zeta \left( 5 \right)\nonumber\\
& + \frac{5}{2}{\zeta ^2}\left( 3 \right)\zeta \left( 4 \right) + 2\zeta \left( 2 \right){S_{2,6}} + 7{S_{2,8}} - 2\zeta \left( {2,6,1,1} \right),
\end{align}
\begin{align}\label{equ:6.6}
{S_{{1^2}3,5}} =&  - \frac{{10471}}{{160}}\zeta \left( {10} \right) + 34{\zeta ^2}\left( 5 \right) + 37\zeta \left( 3 \right)\zeta \left( 7 \right) + 5\zeta \left( 2 \right)\zeta \left( 3 \right)\zeta \left( 5 \right)\nonumber\\
& - \frac{5}{2}{\zeta ^2}\left( 3 \right)\zeta \left( 4 \right) - \frac{9}{2}\zeta \left( 2 \right){S_{2,6}} - \frac{{27}}{2}{S_{2,8}} - 2\zeta \left( {3,5,1,1} \right),
\end{align}
\begin{align}\label{equ:6.7}
{S_{{1^2}4,4}} = &\frac{{1028}}{{15}}\zeta \left( {10} \right) - 38{\zeta ^2}\left( 5 \right) - 50\zeta \left( 3 \right)\zeta \left( 7 \right) + 10\zeta \left( 2 \right)\zeta \left( 3 \right)\zeta \left( 5 \right)\nonumber\\
 &- 2{\zeta ^2}\left( 3 \right)\zeta \left( 4 \right) - 5\zeta \left( 2 \right){S_{2,6}} + 25{S_{2,8}} - 2\zeta \left( {4,4,1,1} \right),
\end{align}
\begin{align}\label{equ:6.8}
{S_{{1^2}5,3}} = & - \frac{{17087}}{{480}}\zeta \left( {10} \right) + \frac{{123}}{4}{\zeta ^2}\left( 5 \right) + \frac{{103}}{2}\zeta \left( 3 \right)\zeta \left( 7 \right) - 26\zeta \left( 2 \right)\zeta \left( 3 \right)\zeta \left( 5 \right)\nonumber\\
& + \frac{{25}}{2}{\zeta ^2}\left( 3 \right)\zeta \left( 4 \right) - \frac{{103}}{4}{S_{2,8}} - 2\zeta \left( {5,3,1,1} \right),
\end{align}
\begin{align}\label{equ:6.9}
{S_{{1^2}6,2}} = & - \frac{{149}}{{40}}\zeta \left( {10} \right) + \frac{9}{2}{\zeta ^2}\left( 5 \right) + 6\zeta \left( 3 \right)\zeta \left( 7 \right) - 3{\zeta ^2}\left( 3 \right)\zeta \left( 4 \right)\nonumber\\
& - 6\zeta \left( 2 \right){S_{2,6}} + 11{S_{2,8}} - 2\zeta \left( {6,2,1,1} \right).
\end{align}
Taking $p=3,m=5$ in Lemma \ref{lem:6.1} and using the duality relation (\ref{equ:6.3}), then
\begin{align*}
&\frac{1}{2}{S_{{1^3},7}} - \frac{1}{2}{S_{12,7}} - {S_{{1^2},8}} + {S_{1,9}}\\
 = &\frac{{2123}}{{160}}\zeta \left( {10} \right) - \frac{{27}}{4}{\zeta ^2}\left( 5 \right) - \frac{{27}}{2}\zeta \left( 3 \right)\zeta \left( 7 \right) + 3\zeta \left( 2 \right)\zeta \left( 3 \right)\zeta \left( 5 \right)\\
 & + \frac{{15}}{8}{\zeta ^2}\left( 3 \right)\zeta \left( 4 \right) + \frac{5}{4}{S_{2,8}}\\
 = &4\zeta \left( {5,{{\left\{ 1 \right\}}_5}} \right) + \sum\limits_{i = 1}^5 {\zeta \left( {4,{{\left\{ 1 \right\}}_{i - 1}}2,{{\left\{ 1 \right\}}_{5 - i}}} \right)} \\
 =& 4\zeta \left( {5,{{\left\{ 1 \right\}}_5}} \right) + \zeta \left( {2,6,1,1} \right) + \zeta \left( {3,5,1,1} \right) + \zeta \left( {4,4,1,1} \right) + \zeta \left( {5,3,1,1} \right) + \zeta \left( {6,2,1,1} \right).
\end{align*}
From \cite{Xu2017} we can get
\[\zeta \left( {5,{{\left\{ 1 \right\}}_5}} \right) = \frac{{973}}{{160}}\zeta \left( {10} \right) - 4{\zeta ^2}\left( 5 \right) - 8\zeta \left( 3 \right)\zeta \left( 7 \right) + 3\zeta \left( 2 \right)\zeta \left( 3 \right)\zeta \left( 5 \right) + \frac{9}{8}{\zeta ^2}\left( 3 \right)\zeta \left( 4 \right).\]
Thus, we obtain
\begin{align}\label{equ:6.10}
&\zeta \left( {2,6,1,1} \right) + \zeta \left( {3,5,1,1} \right) + \zeta \left( {4,4,1,1} \right) + \zeta \left( {5,3,1,1} \right) + \zeta \left( {6,2,1,1} \right)\nonumber\\
= & - \frac{{1769}}{{160}}\zeta \left( {10} \right) + \frac{{37}}{4}{\zeta ^2}\left( 5 \right) + \frac{{27}}{2}\zeta \left( 3 \right)\zeta \left( 7 \right) - 9\zeta \left( 2 \right)\zeta \left( 3 \right)\zeta \left( 5 \right)\nonumber\\
& - \frac{{21}}{8}{\zeta ^2}\left( 3 \right)\zeta \left( 4 \right) + \frac{5}{4}{S_{2,8}}.
\end{align}
In \cite{EW2012}, Eie and Wei proved that for positive integer $p>1$, all mzvs $\zeta \left( {p,p,1,1} \right)$ can be expressed in terms of mzvs of depth $\leq 3$, and gave explicit formulas. Hence, we obtain
\begin{align}\label{equ:6.11}
\zeta \left( {4,4,1,1} \right) =& \frac{{9649}}{{960}}\zeta \left( {10} \right) - 11{\zeta ^2}\left( 5 \right) - \frac{{275}}{{16}}\zeta \left( 3 \right)\zeta \left( 7 \right) + 12\zeta \left( 2 \right)\zeta \left( 3 \right)\zeta \left( 5 \right)\nonumber\\
& - \frac{7}{4}{\zeta ^2}\left( 3 \right)\zeta \left( 4 \right) - \frac{{15}}{4}\zeta \left( 2 \right){S_{2,6}} + 7{S_{2,8}}.
\end{align}
Thus, substituting (\ref{equ:6.11}) into (\ref{equ:6.7}), we can prove the identity (\ref{equ:15}).

It is now proven that all MZVs of weight up to 12 are reducible to {\bf Q}-linear combinations of $\zeta(5,3),\zeta(7,3),\zeta(3,5,3),\zeta(9,3),\zeta(\bar 7,\bar 5)$, single zeta values, and products of these terms, see \cite{BBV2010,Bro2013}. For example, Jos Vermaseren derived this proven identity from {\bf DataMine}:
\begin{align*}
\zeta \left( {6,4,1,1} \right) =&  - \frac{{64}}{{27}}\zeta \left( {7,5} \right) - \frac{{64}}{{27}}\zeta \left( {\bar 7,\bar 5} \right) - \frac{{7967}}{{1944}}\zeta \left( {9,3} \right) + \frac{1}{{12}}{\zeta ^4}\left( 3 \right)\\
& + \frac{{11431}}{{1296}}\zeta \left( 5 \right)\zeta \left( 7 \right) - \frac{{799}}{{72}}\zeta \left( 3 \right)\zeta \left( 9 \right) + 3\zeta \left( 2 \right)\zeta \left( {7,3} \right)\\
& + \frac{7}{2}\zeta \left( 2 \right){\zeta ^2}\left( 5 \right) + 10\zeta \left( 2 \right)\zeta \left( 3 \right)\zeta \left( 7 \right) + \frac{3}{2}\zeta \left( 4 \right)\zeta \left( {5,3} \right)\\
& - \frac{1}{2}\zeta \left( 3 \right)\zeta \left( 4 \right)\zeta \left( 5 \right) - \frac{9}{4}{\zeta ^2}\left( 3 \right)\zeta \left( 6 \right) - \frac{{196274855}}{{10746432}}\zeta \left( {12} \right).
\end{align*}
J. Zhao (personal communication) has obtained the closed forms of MZVs $\zeta \left( {2,6,1,1} \right),\zeta \left( {3,5,1,1} \right)$ and $\zeta \left( {5,3,1,1} \right)$ by using the method of experimental evaluation of {\bf MAPLE}:
\begin{align}\label{equ:6.12}
\zeta \left( {2,6,1,1} \right) =& \frac{{18719}}{{640}}\zeta \left( {10} \right) - \frac{{459}}{{32}}{\zeta ^2}\left( 5 \right) - \frac{{697}}{{32}}\zeta \left( 3 \right)\zeta \left( 7 \right) + \zeta \left( 2 \right)\zeta \left( 3 \right)\zeta \left( 5 \right)\nonumber\\
& + \frac{{23}}{8}{\zeta ^2}\left( 3 \right)\zeta \left( 4 \right) + \zeta \left( 2 \right){S_{2,6}} + \frac{{275}}{{64}}{S_{2,8}},
\end{align}
\begin{align}\label{equ:6.13}
\zeta \left( {3,5,1,1} \right) =&  - \frac{{2671}}{{80}}\zeta \left( {10} \right) + \frac{{315}}{{16}}{\zeta ^2}\left( 5 \right) + \frac{{393}}{{16}}\zeta \left( 3 \right)\zeta \left( 7 \right) - 4\zeta \left( 2 \right)\zeta \left( 3 \right)\zeta \left( 5 \right)\nonumber\\
&- \frac{9}{8}{\zeta ^2}\left( 3 \right)\zeta \left( 4 \right) - \frac{7}{4}\zeta \left( 2 \right){S_{2,6}} - \frac{{147}}{{32}}{S_{2,8}},
\end{align}
\begin{align}\label{equ:6.14}
\zeta \left( {5,3,1,1} \right) = & - \frac{{58}}{{15}}\zeta \left( {10} \right) + \frac{{113}}{{16}}{\zeta ^2}\left( 5 \right) + \frac{{341}}{{16}}\zeta \left( 3 \right)\zeta \left( 7 \right) - \frac{{31}}{2}\zeta \left( 2 \right)\zeta \left( 3 \right)\zeta \left( 5 \right)\nonumber\\
& - \frac{3}{8}{\zeta ^2}\left( 3 \right)\zeta \left( 4 \right) + \frac{{25}}{4}\zeta \left( 2 \right){S_{2,6}} - \frac{{245}}{{32}}{S_{2,8}}.
\end{align}
From (\ref{equ:6.10}), (\ref{equ:6.11}) and (\ref{equ:6.12})-(\ref{equ:6.14}), we also obtain the result
\begin{align}\label{equ:6.15}
\zeta \left( {6,2,1,1} \right) = & - \frac{{1677}}{{128}}\zeta \left( {10} \right) + \frac{{251}}{{32}}{\zeta ^2}\left( 5 \right) + \frac{{371}}{{32}}\zeta \left( 3 \right)\zeta \left( 7 \right) - \frac{5}{2}\zeta \left( 2 \right)\zeta \left( 3 \right)\zeta \left( 5 \right)\nonumber\\
& - \frac{9}{4}{\zeta ^2}\left( 3 \right)\zeta \left( 4 \right) - \frac{7}{4}\zeta \left( 2 \right){S_{2,6}} + \frac{{141}}{{64}}{S_{2,8}}.
\end{align}
\begin{table}[htbp]
 \begin{tabular}{lll}
  \hline
  MZVs & Numerical values of closed form & Approx. value of mzvs\\
 \hline
 $\zeta(2,6,1,1)$ & 0.00027686227249488302455893961&0.00027686227249488302455894055\\
 $\zeta(3,5,1,1)$  &0.00011903399855304478279509062 &0.00011903399855304478279508911  \\
 $\zeta(4,4,1,1)$  &0.00006792914736628345920749066 &0.00006792914736628345920748814  \\
 $\zeta(5,3,1,1)$  &0.00004342997243903687018023324 &0.00004342997243903687018023761\\
$\zeta(6,2,1,1)$  &0.00002950604505239677926788846 &0.00002950604505239677926788724\\
  \hline
 \end{tabular}
 \begin{center}
  \textbf\ \ {\bf TABLE 1.} Numerical values of MZVs
  \end{center}
\end{table}
Based on the results of Table 1, we believe that the obtained closed form representations of MZVs are true.
Hence, we may deduce from equations (\ref{equ:6.5})-(\ref{equ:6.9}) and (\ref{equ:6.12})-(\ref{equ:6.15}) the desired results (\ref{equ:13})-(\ref{equ:17}). Furthermore, letting $m=3$ in formula (2.41) of \cite{X2016}, we find that
\begin{align}\label{equ:6.16}
{S_{{1^4},6}} - 6{S_{{1^2}2,6}} =& \sum\limits_{n = 1}^\infty  {\frac{{H_n^4 - 6H_n^2H_n^{\left( 2 \right)}}}{{{n^6}}}}\nonumber\\
=&\frac{{37949}}{{160}}\zeta \left( {10} \right) - 113{\zeta ^2}\left( 5 \right) - 225\zeta \left( 3 \right)\zeta \left( 7 \right) + 24\zeta \left( 2 \right)\zeta \left( 3 \right)\zeta \left( 5 \right)\nonumber\\
& + 45{\zeta ^2}\left( 3 \right)\zeta \left( 4 \right) + 2\zeta \left( 2 \right){S_{2,6}} + \frac{{57}}{2}{S_{2,8}}.
\end{align}
Substituting (\ref{equ:13}) into (\ref{equ:6.16}) we have the result (\ref{equ:12}).

\section{Further results}

In fact, by using the methods of this paper and reference \cite{X2016,X2017,Xu2017}, it is possible to evaluate other nonlinear Euler sums involving harmonic numbers. For example, from Theorem 3.7, we can establish the relation
\begin{align}\label{equ:6.17}
{S_{{2^2}3,3}} + {S_{{3^2}2,2}} =& \sum\limits_{n = 1}^\infty  {\frac{{{{\left( {H_n^{\left( 2 \right)}} \right)}^2}H_n^{\left( 3 \right)}}}{{{n^3}}} + \sum\limits_{n = 1}^\infty  {\frac{{{{\left( {H_n^{\left( 3 \right)}} \right)}^2}H_n^{\left( 2 \right)}}}{{{n^2}}}} }\nonumber \\
 = & - \frac{{{\rm{1991}}}}{{40}}\zeta \left( {10} \right) + \frac{{83}}{2}{\zeta ^2}\left( 5 \right) + 77\zeta \left( 3 \right)\zeta \left( 7 \right) - 42\zeta \left( 2 \right)\zeta \left( 3 \right)\zeta \left( 5 \right)\nonumber\\
  &+ \frac{{13}}{4}{\zeta ^2}\left( 3 \right)\zeta \left( 4 \right) + 6\zeta \left( 2 \right){S_{2,6}} - \frac{{27}}{2}{S_{2,8}}.
\end{align}
From (2.24) of \cite{X2016}, (2.22) of \cite{X2017} and (3.23) of \cite{Xu2017}, we can confirm that the following three nonlinear Euler sums of weight equal to 10
\[{S_{{1^5},5}} = \sum\limits_{n = 1}^\infty  {\frac{{H_n^5}}{{{n^5}}}},\ {S_{{1^3}2,5}} = \sum\limits_{n = 1}^\infty  {\frac{{H_n^3H_n^{\left( 2 \right)}}}{{{n^5}}}} \quad {\rm and}\quad {S_{{{12}^2},5}} = \sum\limits_{n = 1}^\infty  {\frac{{{H_n}{{\left( {H_n^{\left( 2 \right)}} \right)}^2}}}{{{n^5}}}} \]
can be expressed as a rational linear combination of products of linear sums and zeta values.

Next, we give the closed form of Euler sums $S_{1^5,5},S_{1^32,5}$ and $S_{12^2,5}$. We need establish three equations involving Euler sums.
Letting $r=4,s=3,t=1$ in
(2.22) of \cite{X2017}, we get
\begin{align}\label{equ:6.18}
 &5{S_{{1^2}2,6}} + 2{S_{{1^2}3,5}} + {S_{{{12}^2},5}} - 10{S_{13,6}} - 5{S_{14,5}}\nonumber \\
  &=  - \frac{9}{2}{S_{{1^3},7}} + \frac{{27}}{2}{S_{12,7}} + 14{S_{{1^2},8}} - 9{S_{1,2}}{S_{1,6}} + {S_{1,3}}{S_{1,5}}\nonumber \\
  &\quad- 5\zeta (2){S_{{1^2},6}} + 2\zeta (3){S_{{1^2},5}} - 2\zeta (2){S_{12,5}}.
\end{align}
Combining related equation, we deduce that
\begin{align}\label{equ:6.19}
 &5{S_{{1^2}2,6}} + 2{S_{{1^2}3,5}} + {S_{{{12}^2},5}} \nonumber\\
  &=  - \frac{{3803}}{{40}}\zeta (10) + \frac{{91}}{2}{\zeta ^2}(5) + 92\zeta (3)\zeta (7) - 3\zeta (2)\zeta (3)\zeta (5)\nonumber \\
 &\quad - 20{\zeta ^2}(3)\zeta (4) - \zeta (2){S_{2,6}} - \frac{{35}}{2}{S_{2,8}}.
\end{align}
Setting $k=5,m=4$ in (2.24) of \cite{X2016} yields
\begin{align}\label{equ:6.20}
3{S_{{{12}^2},5}} - 2{S_{{1^3}2,5}} =& \frac{{15923}}
{{80}}\zeta \left( {10} \right) - \frac{{1247}}
{{16}}{\zeta ^2}\left( 5 \right) - \frac{{2293}}
{{16}}\zeta \left( 3 \right)\zeta \left( 7 \right) - 27\zeta \left( 2 \right)\zeta \left( 3 \right)\zeta \left( 5 \right)\nonumber\\
& + \frac{{177}}
{4}{\zeta ^2}\left( 3 \right)\zeta \left( 4 \right) + 8\zeta \left( 2 \right){S_{2,6}} + \frac{{1015}}
{{32}}{S_{2,8}}.
\end{align}
In \cite{KO2010}, Masanobu Kaneko and Yasuo Ohno proved that
\begin{align}\label{equ:6.21}
&{\left( { - 1} \right)^k}{\zeta ^ \star }\left( {k + 1,{{\left\{ 1 \right\}}_n}} \right) - {\left( { - 1} \right)^n}{\zeta ^ \star }\left( {n + 1,{{\left\{ 1 \right\}}_k}} \right)\nonumber\\
& = k\zeta \left( {k + 2,{{\left\{ 1 \right\}}_{n - 1}}} \right) - n\zeta \left( {n + 2,{{\left\{ 1 \right\}}_{k - 1}}} \right)\nonumber\\
&\quad + {\left( { - 1} \right)^k}\sum\limits_{j = 0}^{k - 2} {{{\left( { - 1} \right)}^j}\zeta \left( {k - j} \right)\zeta \left( {n + 1,{{\left\{ 1 \right\}}_j}} \right)}\nonumber \\
&\quad - {\left( { - 1} \right)^n}\sum\limits_{j = 0}^{n - 2} {{{\left( { - 1} \right)}^j}\zeta \left( {n - j} \right)\zeta \left( {k + 1,{{\left\{ 1 \right\}}_j}} \right)} .
\end{align}
Putting $k=4,n=5$ in (\ref{equ:6.21}) we obtain
\begin{align}\label{equ:6.22}
{\zeta ^ \star }\left( {5,{{\left\{ 1 \right\}}_5}} \right) + {\zeta ^ \star }\left( {6,{{\left\{ 1 \right\}}_4}} \right) = &4\zeta \left( {6,{{\left\{ 1 \right\}}_4}} \right) - 5\zeta \left( {7,{{\left\{ 1 \right\}}_3}} \right) + \zeta \left( 4 \right)\zeta \left( 6 \right) - \zeta \left( 3 \right)\zeta \left( {6,1} \right)\nonumber\\
& + \zeta \left( 2 \right)\zeta \left( {6,1,1} \right) + {\zeta ^2}\left( 5 \right) - \zeta \left( 4 \right)\zeta \left( {5,1} \right) 
\nonumber\\&+ \zeta \left( 3 \right)\zeta \left( {5,1,1} \right) - \zeta \left( 2 \right)\zeta \left( {5,{{\left\{ 1 \right\}}_3}} \right).
\end{align}
From \cite{Xu2017}, we know that for any $m,k\in \N$, the multiple zeta value $\zeta(m+1,\{1\}_{k-1} )$ can be
represented as a polynomial of zeta values with rational coefficients. For example:
\[\begin{array}{l}
 \zeta \left( {2,{{\left\{ 1 \right\}}_m}} \right) = \zeta \left( {m + 2} \right), \\
 \zeta \left( {3,{{\left\{ 1 \right\}}_m}} \right) = \frac{{m + 2}}{2}\zeta \left( {m + 3} \right) - \frac{1}{2}\sum\limits_{k = 1}^m {\zeta \left( {k + 1} \right)\zeta \left( {m + 2 - k} \right)} . \\
 \end{array}\]
Moreover, we note that the mzv ${\zeta ^ \star }\left( {6,{{\left\{ 1 \right\}}_4}} \right)$ is equal to 
\begin{align}\label{equ:6.23}
{\zeta ^ \star }\left( {6,{{\left\{ 1 \right\}}_4}} \right) &= \sum\limits_{n = 1}^\infty  {\frac{{\zeta _n^ \star \left( {{{\left\{ 1 \right\}}_4}} \right)}}
{{{n^6}}}}\nonumber \\
& = \frac{1}
{{24}}\sum\limits_{n = 1}^\infty  {\frac{{H_n^4 + 8{H_n}H_n^{\left( 3 \right)} + 6H_n^2H_n^{\left( 2 \right)} + 3{{\left( {H_n^{\left( 2 \right)}} \right)}^2} + 6H_n^{\left( 4 \right)}}}
{{{n^6}}}}\nonumber \\
& = \frac{{14221}}
{{1280}}\zeta \left( {10} \right) - \frac{{221}}
{{32}}{\zeta ^2}\left( 5 \right) - \frac{{343}}
{{32}}\zeta \left( 3 \right)\zeta \left( 7 \right) + 4\zeta \left( 2 \right)\zeta \left( 3 \right)\zeta \left( 5 \right)\nonumber\\&\quad + \frac{3}
{4}{\zeta ^2}\left( 3 \right)\zeta \left( 4 \right) + \frac{{61}}
{{64}}{S_{2,8}}.
\end{align}
Hence, from (\ref{equ:6.22}) and (\ref{equ:6.23}), it is easily show that
\begin{align}\label{equ:6.24}
{\zeta ^ \star }\left( {5,{{\left\{ 1 \right\}}_5}} \right) = & - \frac{{7681}}
{{640}}\zeta \left( {10} \right) + \frac{{253}}
{{32}}{\zeta ^2}\left( 5 \right) + \frac{{407}}
{{32}}\zeta \left( 3 \right)\zeta \left( 7 \right)\nonumber \\&- 3\zeta \left( 2 \right)\zeta \left( 3 \right)\zeta \left( 5 \right) - \frac{{19}}
{8}{\zeta ^2}\left( 3 \right)\zeta \left( 4 \right) - \frac{{61}}
{{64}}{S_{2,8}}.
\end{align}
Similarly, the mzv ${\zeta ^ \star }\left( {5,{{\left\{ 1 \right\}}_5}} \right)$ can be rewritten as
\begin{align}\label{equ:6.25} 
&{\zeta ^ \star }\left( {5,{{\left\{ 1 \right\}}_5}} \right) = \sum\limits_{n = 1}^\infty  {\frac{{\zeta _n^ \star \left( {{{\left\{ 1 \right\}}_5}} \right)}}
{{{n^5}}}}\nonumber \\
& = \frac{1}
{{120}}\sum\limits_{n = 1}^\infty  {\frac{{H_n^5 + 10H_n^3H_n^{\left( 2 \right)} + 20H_n^2H_n^{\left( 3 \right)} + 15{H_n}{{\left( {H_n^{\left( 2 \right)}} \right)}^2} + 30{H_n}H_n^{\left( 4 \right)} + 20H_n^{\left( 2 \right)}H_n^{\left( 3 \right)} + 24H_n^{\left( 5 \right)}}}
{{{n^5}}}} .
\end{align}
So, we obtain the equation
\begin{align}\label{equ:6.26} 
3{S_{{{12}^2},5}} + 2{S_{{1^3}2,5}} = & - \frac{{49069}}
{{160}}\zeta \left( {10} \right) + \frac{{2237}}
{{16}}{\zeta ^2}\left( 5 \right) + \frac{{3103}}
{{16}}\zeta \left( 3 \right)\zeta \left( 7 \right) + 33\zeta \left( 2 \right)\zeta \left( 3 \right)\zeta \left( 5 \right)\nonumber\\
&- \frac{{255}}
{4}{\zeta ^2}\left( 3 \right)\zeta \left( 4 \right) - 2\zeta \left( 2 \right){S_{2,6}} - \frac{{1189}}
{{32}}{S_{2,8}}.
\end{align}
Hence, combining the equations (\ref{equ:6.19}), (\ref{equ:6.20}) and (\ref{equ:6.26}), by a direct calculation, the results are as follows
\begin{align*}
{S_{{1^5},5}} =& \frac{{10089}}{{32}}\zeta \left( {10} \right) - \frac{{1997}}{{16}}{\zeta ^2}\left( 5 \right) - \frac{{3215}}{{16}}\zeta \left( 3 \right)\zeta \left( 7 \right) - 35\zeta \left( 2 \right)\zeta \left( 3 \right)\zeta \left( 5 \right)\\
& + \frac{{215}}{4}{\zeta ^2}\left( 3 \right)\zeta \left( 4 \right) + 5\zeta \left( 2 \right){S_{2,6}} + \frac{{1365}}{{32}}{S_{2,8}},\\
 {S_{{1^3}2,5}} = & - \frac{{16183}}{{128}}\zeta \left( {10} \right) + \frac{{871}}{{16}}{\zeta ^2}\left( 5 \right) + \frac{{1349}}{{16}}\zeta \left( 3 \right)\zeta \left( 7 \right) + 15\zeta \left( 2 \right)\zeta \left( 3 \right)\zeta \left( 5 \right) \\
  &- 27{\zeta ^2}\left( 3 \right)\zeta \left( 4 \right) - \frac{5}{2}\zeta \left( 2 \right){S_{2,6}} - \frac{{551}}{{32}}{S_{2,8}}, \\
 {S_{{{12}^2},5}} = & - \frac{{5741}}{{320}}\zeta \left( {10} \right) + \frac{{165}}{{16}}{\zeta ^2}\left( 5 \right) + \frac{{135}}{{16}}\zeta \left( 3 \right)\zeta \left( 7 \right) + \zeta \left( 2 \right)\zeta \left( 3 \right)\zeta \left( 5 \right) \\
  &- \frac{{13}}{4}{\zeta ^2}\left( 3 \right)\zeta \left( 4 \right) + \zeta \left( 2 \right){S_{2,6}} - \frac{{29}}{{32}}{S_{2,8}}.
\end{align*}
In \cite{Xu-2017}, we establish the relation
\begin{align}\label{equ:7.11} 
{S_{{{12}^2},5}} + 2{S_{125,2}} &=  - \frac{{2403}}{{160}}\zeta \left( {10} \right) + \frac{{69}}{4}{\zeta ^2}\left( 5 \right) + \frac{{491}}{8}\zeta \left( 3 \right)\zeta \left( 7 \right) - 34\zeta \left( 2 \right)\zeta \left( 3 \right)\zeta \left( 5 \right)\nonumber\\
&\quad- \frac{{27}}{4}{\zeta ^2}\left( 3 \right)\zeta \left( 4 \right) + \frac{{27}}{2}\zeta \left( 2 \right){S_{2,6}} - 13{S_{2,8}},
\end{align}
Hence, we obtain the closed form
\begin{align*}
 {S_{125,2}} &= \frac{{187}}{{128}}\zeta \left( {10} \right) + \frac{{111}}{{32}}{\zeta ^2}\left( 5 \right) + \frac{{847}}{{32}}\zeta \left( 3 \right)\zeta \left( 7 \right) - \frac{{35}}{2}\zeta \left( 2 \right)\zeta \left( 3 \right)\zeta \left( 5 \right) \\
  &\quad- \frac{7}{4}{\zeta ^2}\left( 3 \right)\zeta \left( 4 \right) + \frac{{25}}{4}\zeta \left( 2 \right){S_{2,6}} - \frac{{387}}{{64}}{S_{2,8}}.
\end{align*}
Furthermore, we can also deduce the following two closed form of Euler sums ${S_{{1^4}2,4}}$ and ${S_{123,4}}$ that
\begin{align*}
 {S_{{1^4}2,4}} = & - \frac{{48071}}{{480}}\zeta \left( {10} \right) + \frac{{569}}{8}{\zeta ^2}\left( 5 \right) + \frac{{429}}{8}\zeta \left( 3 \right)\zeta \left( 7 \right) + 28\zeta \left( 2 \right)\zeta \left( 3 \right)\zeta \left( 5 \right) \\
  &- \frac{{97}}{2}{\zeta ^2}\left( 3 \right)\zeta \left( 4 \right) - 4\zeta \left( 2 \right){S_{2,6}} - \frac{{225}}{{16}}{S_{2,8}}, \\
{S_{123,4}} =&  - \frac{{2849}}{{192}}\zeta \left( {10} \right) + \frac{{59}}{8}{\zeta ^2}\left( 5 \right) + \frac{{303}}{{16}}\zeta \left( 3 \right)\zeta \left( 7 \right) - 5\zeta \left( 2 \right)\zeta \left( 3 \right)\zeta \left( 5 \right) \\
  &+ {\zeta ^2}\left( 3 \right)\zeta \left( 4 \right) - \frac{7}{4}\zeta \left( 2 \right){S_{2,6}} - \frac{{51}}{{16}}{S_{2,8}}.
\end{align*}
Moreover, we use Mathematica tool to check numerically each of the specific identities listed. The numerical values of
nonlinear Euler sums of weights 10, to 30 decimal digits, see Table 2.

Finally, we close this paper with two conjectures.
\begin{con}\label{con1}  All non-alternating Euler sums of weight ten can be expressed as a rational linear combination of $\zeta \left( {10} \right),{\zeta ^2}\left( 5 \right),\zeta \left( 3 \right)\zeta \left( 7 \right),\zeta \left( 2 \right)\zeta \left( 3 \right)\zeta \left( 5 \right),{\zeta ^2}\left( 3 \right)\zeta \left( 4 \right),\zeta \left( 2 \right){S_{2,6}}$ and ${S_{2,8}}$.
\end{con}

\begin{con}\label{con2} For positive integer $m$, the combined sum 
\[\zeta \left( {{\overline {2m}},2} \right) + 2m\zeta \left( {{\overline {2m + 1}},1} \right)\]
 can be expressed in terms of non-alternating double zeta values.
\end{con}
In this paper, we can prove that if $m=1,2,3$, the Conjecture \ref{con2} are correct. We obtain the following three closed form evaluations
\begin{align*}
 &\zeta \left( {\bar 2,2} \right) + 2\zeta \left( {\bar 3,1} \right) = \frac{5}{{16}}\zeta \left( 4 \right), \\
 &\zeta \left( {\bar 4,2} \right) + 4\zeta \left( {\bar 5,1} \right) = \frac{1}{4}\zeta \left( {4,2} \right) - \frac{{37}}{{64}}\zeta \left( 6 \right) + \frac{1}{2}{\zeta ^2}\left( 3 \right), \\
 &\zeta \left( {\bar 6,2} \right) + 6\zeta \left( {\bar 7,1} \right) =  - \frac{9}{8}\zeta \left( {6,2} \right) - \frac{{1887}}{{256}}\zeta \left( 8 \right) + 6\zeta \left( 3 \right)\zeta \left( 5 \right).
\end{align*}
The experimental evaluation (\ref{equ:5.5}) further confirm that the Conjecture \ref{con2} is correct when $m=4$.

{\bf Acknowledgments.} The author would like to thank  professor J. Zhao for sending us the experimental evaluations (\ref{equ:6.12})-(\ref{equ:6.14}) upon which we can build our current work.

 \newcommand{\tabincell}[2]{\begin{tabular}{@{}#1@{}}#2\end{tabular}}
\begin{table}[htbp]
 \begin{tabular}{|c|c|c|}
  \hline
 \tabincell{c}{Euler sum} &\tabincell{c}{ Numerical values of closed form\\ (30 decimal digits)} & \tabincell{c}{Numerical approximation of\\ Euler sum (30 decimal digits)}\\
  \hline
$S_{26,2}$&1.90896974096899176135265551253&1.90896974096899176135265550729\\
\hline
$S_{23,5}$&1.05312633310471200888384448874&1.05312633310471200888384448505\\
\hline
$S_{25,3}$&1.27456762963066440047814901707&1.27456762963066440047814902209\\
\hline
$S_{35,2}$&1.77427138402017597709644693181&1.77427138402017597709644692042\\
\hline
$S_{24,4}$&1.11283733496948834629784104857&1.11283733496948834629784103882\\
\hline
$S_{4^2,2}$&1.74226189661777917404814402680&1.74226189661777917404814403681\\
\hline
$S_{3^2,4}$&1.10624811674322343965670197651&1.10624811674322343965670198881\\
\hline
$S_{34,3}$&1.24696428137564267341799534681&1.24696428137564267341799533524\\
\hline
$S_{2^3,4}$&1.17640873941393521751500234937&1.17640873941393521751500235018\\
\hline
$S_{2^4,2}$&3.51103090093403934158187558633&3.51103090093403934158187559047\\
\hline
$S_{2^24,2}$&2.34320172841378735744027226238&2.34320172841378735744027225178\\
\hline
$S_{1^4,6}$& 1.10265621583115731907315836205 &1.10265621583115731907315836163\\
\hline
$S_{1^22,6}$& 1.05260005577214971264422620961 &1.05260005577214971264422620815\\
\hline
$S_{1^23,5}$& 1.10432986872601288241896330567&1.10432986872601288241896330207\\
\hline
$S_{1^24,4}$& 1.23696110536553010766655808317&1.23696110536553010766655808129\\
\hline
$S_{1^25,3}$&1.67416857092838554996910378694 &1.67416857092838554996910378192\\
\hline
$S_{1^26,2}$&4.66117013404954343340003376622&4.66117013404954343340003376276\\
\hline
$S_{1^5,5}$&1.41632664251398820119526325885&1.41632664251398820119526325895\\
\hline
$S_{1^32,5}$&1.19189395441095256742841205632&1.19189395441095256742841204947\\
\hline
$S_{12^2,5}$&1.09450138606515480472596873607&1.09450138606515480472596873429\\
\hline
$S_{125,2}$& 3.08534687169631240730809129046&3.08534687169631240730809129197\\
\hline
$S_{1^42,4}$&1.92786760234663387325396559374&1.92786760234663387325396558568\\
\hline
$S_{123,4}$&1.19726207517224197717128588358&1.19726207517224197717128587807\\
\hline
 \end{tabular}\\
 [2mm]
  \centering
{\bf TABLE 2.} Numerical approximation of Euler sums of weight=10
\end{table}
\end{document}